\documentclass[12pt]{article}

\usepackage{latexsym}
\usepackage{calc}
\usepackage{datetime}
\usepackage{amssymb, amsmath, amsfonts, listings, mathrsfs} 
\usepackage{xcolor} 
\usepackage[normalem]{ulem} 

\usepackage{graphicx}
\usepackage[labelfont=bf]{caption}

\newcounter{hours}\newcounter{minutes}

\newcommand{\stkout}[1]{\ifmmode\text{\sout{\ensuremath{#1}}}\else\sout{#1}\fi} 

\def\nr{\par \noindent}

\def\Def{\stackrel{\mathrm{def}}{=}}

\def\inter{{\rm int \,}}

\def\dom{{\rm dom \,}}
\def\vf{\varphi}
\def\beq{\begin{equation}}
\def\eeq{\end{equation}}

\def\R{\mathbb{R}}
\def\E{\mathbb{E}}

\def\BI{\begin{itemize}}
\def\EI{\end{itemize}}

\newcommand{\SetEQ}{\setcounter{equation}{0}}
\newcommand{\refLE}[1]{\ensuremath{\stackrel{(\ref{#1})}{\leq}}}
\newcommand{\refEQ}[1]{\ensuremath{\stackrel{(\ref{#1})}{=}}}
\newcommand{\refGE}[1]{\ensuremath{\stackrel{(\ref{#1})}{\geq}}}

\newtheorem{theorem}{Theorem}
\newtheorem{lemma}{Lemma}
\newtheorem{corollary}{Corollary}

\newtheorem{assumption}{Assumption}
\newtheorem{definition}{Definition}

\newtheorem{example}{Example}
\newtheorem{remark}{Remark}
\newcommand{\proof}{\bf Proof: \rm \nr}
\newcommand{\qed}{\hfill $\Box$ \nr \medskip}
\newcommand{\half}{\mbox{${1 \over 2}$}}

\def\ba{\begin{array}}
\def\ea{\end{array}}
\def\beann{\begin{eqnarray*}}
\def\eeann{\end{eqnarray*}}
\def\bea{\begin{eqnarray}}
\def\eea{\end{eqnarray}}

\def\BT{\begin{theorem}}
\def\ET{\end{theorem}}
\def\BL{\begin{lemma}}
\def\EL{\end{lemma}}
\def\BC{\begin{corollary}}
\def\EC{\end{corollary}}
\def\BE{\begin{example}}
\def\EE{\end{example}}
\def\BD{\begin{definition}}
\def\ED{\end{definition}}
\def\BR{\begin{remark}}
\def\ER{\end{remark}}
\def\BAS{\begin{assumption}}
\def\EAS{\end{assumption}}
\def\BI{\begin{itemize}}
\def\EI{\end{itemize}}

\def\BMP{\begin{minipage}{9.5cm}}
\def\EMP{\end{minipage}}
\def\MPT{\begin{minipage}{11.5cm}}
\def\EPT{\end{minipage}}

\def\la{\langle}
\def\ra{\rangle}

\def\QF{\hspace{5ex} \Box}
\def\QR{\hfill \Box}

\title{Primal subgradient methods with predefined stepsizes}
\author{Yu. Nesterov
\thanks{Center for Operations Research and Econometrics (CORE),
Catholic University of Louvain (UCL). E-mail:
Yurii.Nesterov@uclouvain.be. This paper has received funding from the European Research Council (ERC) under the European Union’s Horizon 2020 research and innovation program (grant agreement No 788368).
It was also supported by Multidisciplinary Institute in Artificial intelligence MIAI@Grenoble Alpes (ANR-19-P3IA-0003).
}}

\begin{document}
\maketitle

\abstract{In this paper, we suggest a new framework for analyzing primal subgradient methods for nonsmooth convex optimization problems. We show that the classical step-size rules, based on normalization of subgradient, or on the knowledge of optimal value of the objective function, need corrections when they are applied to optimization problems with constraints. Their proper modifications allow a significant acceleration of these schemes when the objective function has favorable properties (smoothness, strong convexity). We show how the new methods can be used for solving optmization problems with functional constraints with possibility to approximate the optimal Lagrange multipliers. One of our primal-dual methods works also for unbounded feasible set.}

\vspace{10ex}\noindent
{\bf Keywords:} Convex optimization, nonsmooth optimization, subgradient methods, constrained problems, optimal Lagrange multipliers.



\section{Introduction}
\setcounter{equation}{0}

\vspace{1ex}\noindent
{\bf Motivation.} The first method for unconstrained minimization of  nonsmooth convex function was proposed in \cite{Shor}. This was a primal subgradient method 
\beq\label{met-SGM}
\ba{rcl}
x_{k+1} & = & x_k - h_k d_k, \quad d_k \; = \; f'(x_k), \quad k \geq 0,
\ea
\eeq
with constant step sizes $h_k \equiv h > 0$, where $f'(x_k)$ is a subgradient of the objective function at point $x_k$. In the next years, there were developed several strategies for choosing the steps (see \cite{BT} for historical remarks and references). Among them, the most important one is the rule of the {\em first-order} divergent series:
\beq\label{eq-Div}
\ba{c}
h_k > 0, \quad h_k \to 0, \quad \sum\limits_{k=0}^{\infty} h_k = \infty,
\ea
\eeq
with the optimal choice $h_k = O(k^{-1/2})$. As a variant, it is possible to use in (\ref{met-SGM}) the normalized directions 
\beq\label{eq-DirN}
\ba{rcl}
d_k & = &  {f'(x_k) \over \| f'(x_k) \|}.
\ea
\eeq

Another alternative for the step sizes is based on the known optimal value $f^*$ \cite{BT}:
\beq\label{eq-Opt}
\ba{rcl}
h_k & = & {f(x_k) - f^* \over \| f'(x_k) \|^2}.
\ea
\eeq
In both cases, the corresponding schemes, as applied to functions with bounded subgradients, have the optimal rate of convergence $O(k^{-1/2})$, established for the best value of the objective function observed during the minimization process \cite{NY}. The presence of simple set constraints was treated just by applying to the minimization sequence an Euclidean projection onto the feasible set.

The next important advancement in this area is related to development of {\em mirror descent method} (MDM) (\cite{NY}, see also \cite{NJ}). In this scheme, the main information is accumulated in the dual space in the form of aggregated subgradients. For defining the next test point, this object is mapped ({\em mirrored}) to the primal space by a special prox-function, related to a general norm. 
Thus, we get an important possibility of describing topology of convex sets by the appropriate norms.

After this discovery, during several decades, the research activity in this field was concentrated on the development of dual schemes. One of the drawbacks of the classical MDM is that the new subgradients are accumulated with the vanishing weights $h_k$. It was corrected in the framework of {\em dual averaging} \cite{DA}, where the aggregation coefficients can be even increasing, and the convergence in the primal space is achieved by applying some vanishing scaling coefficients. 
Another drawback is related to the fact that the convergence guarantees are traditionally established only for the best values of the objective function. This inconvenience was eliminated by development of {\em quasi-monotone dual methods} \cite{NS}, where the rate of convergence is proved for all points of the minimization sequence.

Thus, at some moment, primal methods were almost forgotten. However, in this paper we are going to show that in some situations the primal schemes are very useful. Moreover, there is still space for improvement of the classical methods. Our optimism is supported by the following observations.

Firstly, from the recent developments in Optimization Theory, it becomes clear that the size of subgradients of objective functions for the problems with simple set constraints must be defined differently. Hence, the usual norms in the rules (\ref{eq-DirN}) and \ref{eq-Opt}) can be replaced by more appropriate objects.

Secondly, for an important class of {\em quasi-convex functions}, linear approximations do not work properly. Hence, for corresponding optimization problems, only the primal schemes can be used.
Finally, as we will see, the proper primal schemes provide us with a very simple and natural possibility for approximating the optimal Lagrange multipliers for problems with functional constraints, eliminating somehow the heavy machinery, which is typical for the methods based on Augmented Lagrangian.

\vspace{1ex}\noindent
{\bf Contents.} In Section \ref{sc-Quasi}, we present a new subgradient method for minimizing quasi-convex function on a simple set. Its justification is based on a new concept of {\em directional proximity measure}, which is a generalization of the old technique initially presented in \cite{Quasi}. In this method, we apply and {\em indirect strategy} for choosing the step size, which needs solution of a simple univariate equation. 
In unconstrained case, this strategy is reduced to the normalization step (\ref{eq-DirN}). The main advantage of the new method is a possibility of automatic acceleration for functions with H\"older-continuous gradients.

In Section \ref{sc-Conv}, we present a method for solving composite minimization problem with max-type objective function. For choosing the step size, we use a proper  generalization of the rule (\ref{eq-Opt}), based on an optimal value of the objective. This method admits a linear rate of convergence for smooth strongly convex functions (see Section \ref{sc-RelSmooth}). Note that a simple example demonstrates that the classical rule does not benefit from strong convexity. Method of Section \ref{sc-Conv} automatically accelerates on functions with H\"older continuous gradient.

In Sections \ref{sc-Cons}, we consider a minimization problem with single max-type constraint containing an additive composite term. For this problem, we apply a swithching strategy, where the steps for the objective function are based on the rule of Section \ref{sc-Quasi}, and for improving the feasibility we use the step size strategy of Section \ref{sc-Conv}. For controlling the step sizes, we suggest a new rule of the {\em second-order divergent series}:
$$
\ba{rcl}
\tau_k & \geq & \tau_{k+1} \; > \; 0, \quad \sum\limits_{k=0}^{\infty} \tau_k^2 \; = \; \infty.
\ea
$$
For the bounded feasible sets, it eliminates an unpleasant logarithmic factors in the convergence rate. The method automatically accelerates for the problem with smooth functional components. It is interesting that the rates of convergence for the objective function and the constraints can be different.

The remaining sections of the paper are devoted to the methods, which can approximate optimal Lagrange multipliers for convex problems with functional inequality constraints. In Section \ref{sc-Lagrange1}, we consider simplest swithcing strategy of this type, where for the steps with the objective function we use the rule of Section \ref{sc-Quasi}, and for the steps with violated constraints we use the rule of Section \ref{sc-Conv}. In the method of Section \ref{sc-Lagrange2}, both steps are based the rule of Section \ref{sc-Quasi}.
In both cases, we obtain the rates of convergence for infeasibility of the generated points, and the upper bound for the duality gap, computed for the simple estimates of the optimal dual multipliers. Such an estimate is formed as a sum of steps at active iterations for each violated constraint divided by the sum of steps at iterations when the objective function was active.

In Section \ref{sc-AcDual}, we provide the theoretical guarantees for our estimates of the optimal Lagrange multipliers in terms of value of the dual function. They depend on the depth of Slater condition of our problem. Finally, in Section \ref{sc-NoBound}, we present a swithching method, which can generate the approximate dual multipliers for problems with unbounded feasible set.

\vspace{1ex}\noindent
{\bf Notation.} 
Denote by $\E$ a finite-dimensional real vector space,
and by $\E^*$ its dual space composed by linear functions
on $\E$. For such a function $s \in \E^*$, denote by
$\la s, x \ra$ its value at $x \in \E$.
For measuring distances in $\E$, we use an arbitrary norm $\| \cdot \|$. The corresponding dual norm is defined in a standard way:
$$
\ba{rcl}
\| g \|_* \; = \; \max\limits_{x \in \E} \Big\{ \; \la g,x \ra, \; \| x \| \leq 1 \; \Big\}, \quad g \in \E^*.
\ea
$$
Sometimes it is convenient to measure distances in $\E$ by Euclidean norm $\| \cdot \|_B$. It is defined by a self-adjoint positive-definite linear operator $B: \E \to \E^*$ in the following way:
\beq\label{def-Euclid}
\ba{rcl}
\| x \|_B & = & \la B x, x \ra^{1/2}, \; x \in \E, \quad \| g \|^*_B \; = \; \la g, B^{-1} g \ra^{1/2}, \; g \in \E^*.
\ea
\eeq
In case of $\E = \R^n$, for $x \in \R^n$, we use notation $\| x \|^2_2 = \sum\limits_{i=1}^n (x^{(i)})^2$.

For a differentiable function $f(\cdot)$ with convex and open domain $\dom f \subseteq \E$, denote by $\nabla f(x) \in \E^*$ its gradient at point $x \in \dom f$. If $f$ is convex, it can be used for defining the {\em Bregman distance} between two points $x, y \in \dom f$:
\beq\label{def-Breg}
\ba{rcl}
\beta_f(x,y) & = & f(y) - f(x) - \la \nabla f(x), y - x \ra \; \geq \; 0.
\ea
\eeq

In this paper, we develop new {\em proximal-gradient methods} based on a predefined {\em prox-function} $d(\cdot)$, which can be restricted to a convex open domain $\dom d \subseteq \E$. This domain always contains the basic feasible set of the corresponding optimization problem. We assume that $d(\cdot)$ is continuously differentiable and {\em strongly convex} on $\dom d$ with parameter one:
\beq\label{def-Strong}
\ba{rcl}
d(y) & \geq & d(x) + \la \nabla d(x), y - x \ra + \half \| y - x \|^2, \quad x, y \in \dom d.
\ea
\eeq
Thus, combining the definition (\ref{def-Breg}) with inequality (\ref{def-Strong}), we get
\beq\label{eq-BGrow}
\ba{rcl}
\beta_d(x,y) & \geq & \half \| x - y \|^2, \quad x, y \in \dom d.
\ea
\eeq

\section{Subgradient method for quasi-convex problems}\label{sc-Quasi}
\SetEQ

Consider the following constrained optimization problem:
\beq\label{prob-Quasi}
\min\limits_{x \in Q} \; f_0(x),
\eeq
where function $f_0(\cdot)$ is closed and quasi-convex on $\dom f_0$ and the set $Q 
\subseteq \dom f_0$ is closed and convex. Denote by $x^*$ one of the optimal solutions of (\ref{prob-Quasi}) and let $f_0^* = f_0(x^*)$. We assume that  
\beq\label{eq-XInter}
\ba{rcl}
x^* & \in & \inter \dom f_0.
\ea
\eeq

Let us assume that at any point $x \in Q$ it is possible to compute a vector $f_0'(x) \in \E^* \setminus \{0\}$, satisfying the following condition: for any $y \in \E$, we have:
\beq\label{def-Grad}
\ba{rcl}
\la f_0'(x), y - x \ra \geq 0  & \Rightarrow & f_0(y) \; \geq \; f_0(x).
\ea
\eeq
(If $y \not\in \dom f_0$, then $f_0(y) \Def + \infty$.) If $f_0(\cdot)$ is differentiable at $x$, then $f'_0(x) = \nabla f_0(x)$.

In order to justify the rate of convergence of our scheme, we need the following characteristic of problem (\ref{prob-Quasi}):
\beq\label{def-Grow}
\ba{rcl}
\mu(r) & \Def & \sup\limits_{x \in \E} \Big\{ f_0(x) - f_0^*: \; \| x - x^* \| \leq r \Big\}, \quad r \geq 0.
\ea
\eeq
If $y \not\in \dom f_0$ and $r = \| y - x^* \|$, then $\mu(r) = + \infty$. In view of assumption (\ref{eq-XInter}), this function is finite at least in some neighborhood of the point $r=0$.

We say that vector $d \in \E$ defines a {\em direction} in $\E$ if $\| d \| = 1$. If $\la f_0'(x), d \ra > 0$, we call it a {\em recession direction} of function $f_0(\cdot)$ at point $x \in Q$. Using such a direction, we can define the {\em directional proximity measure} of point $x \in Q$ as follows:
\beq\label{def-PM}
\ba{rcl}
\delta_d(x) & \Def & {\la f_0'(x), x - x^* \ra \over \la f_0'(x), d \ra } \; \geq \; 0.
\ea
\eeq
\BL\label{lm-PM}
Let $d$ be a recession direction at point $x \in Q$. Then
\beq\label{eq-PM}
\ba{rcl}
f_0(x) - f_0^* & \leq & \mu(\delta_d(x)).
\ea
\eeq
\EL
\proof
Indeed, let us define $y = x^* + \delta_d(x) d$. Then
$$
\ba{rcl}
\la f_0'(x), y \ra & = & \la f_0'(x), x^* \ra + \delta_d(x) \la f_0'(x), d \ra \; \refEQ{def-PM} \la f_0'(x), x \ra.
\ea
$$
Since $f_0(\cdot)$ is quasi-convex, this means that $f_0(y) \refGE{def-Grad} f_0(x)$. Therefore,
$$
\ba{rcl}
f_0(x) - f_0^* & \leq & f_0(y) - f_0^* \; \refLE{def-Grow} \; \mu (\| y - x^* \|) \; = \; \mu(\delta_d(x)). \QF
\ea
$$

In our analysis, we use the following univariate functions:
\beq\label{def-Phi}
\ba{rcl}
\vf_{\bar x}(\lambda) & = & \max\limits_{x \in Q} \Big\{ \lambda \la f_0'(\bar x), \bar x - x \ra - \beta_d(\bar x, x) \Big\}, \quad \lambda \geq 0,
\ea
\eeq
where $\bar x \in Q$. The optimal solution of this optimization problem is denoted by $T_{\bar x}(\lambda)$. Note that function $\vf_{\bar x}(\cdot)$ is convex and continuously differentiable with the following derivative:
\beq\label{eq-DerPhi}
\ba{rcl}
\vf_{\bar x}'(\lambda) & = & \la f_0'(\bar x), \bar x - T_{\bar x}(\lambda) \ra.
\ea
\eeq
Thus, $\vf_{\bar x}(0) = 0$, $T_{\bar x}(0) = \bar x$, and $\vf'_{\bar x}(0) = 0$. Consequently, $\vf_{\bar x}(\lambda) \geq 0$ for all $\lambda \geq 0$.

The basic subgradient method for solving the problem (\ref{prob-Quasi}) looks as follows.
\beq\label{met-SGM}
\ba{|c|}
\hline \\
\mbox{\bf Basic Subgradient Method}\\
\\
\hline \\
\ba{l}
\mbox{{\bf Initialization.}  Choose $x_0 \in Q$.}\\
\\
\mbox{{\bf $k$th iteration ($k \geq 0$).} a) Choose step-size parameter $\lambda_k > 0$.}\\
\\
\mbox{b) Compute $x_{k+1} = T_{x_k}(\lambda_k) \Def \arg\min\limits_{x \in Q} \Big\{ \lambda_k \la f'_0(x_k), x \ra + \beta(x_k,x) \Big\}$.}\\
\\
\ea\\
\hline
\ea
\eeq

\BL\label{lm-Step}
At each iteration $k \geq 0$ of method (\ref{met-SGM}), for any $x \in Q$, we have
\beq\label{eq-Step}
\ba{rcl}
\beta_d(x_{k+1},x) & \leq & \beta_d(x_k,x ) - \lambda_k \la f_0'(x_k), x_k - x \ra + \vf_{x_k}(\lambda_k).
\ea
\eeq
\EL
\proof
Note that the first-order optimality condition for the problem at Step b) can be written as follows:
\beq\label{eq-COpt}
\ba{rcl}
\la \lambda_k f_0'(x_k) + \nabla d(x_{k+1}) - \nabla d(x_k), x - x_{k+1} \ra & \geq & 0, \quad \forall x \in Q.
\ea
\eeq
Therefore, we have
$$
\ba{rl}
& \beta_d(x_{k+1},x) \; \refEQ{def-Breg} \; d(x) - d(x_{k+1}) - \la \nabla d(x_{k+1}), x - x_{k+1} \ra \\
\\
= & \beta_d(x_{k},x) + d(x_{k}) + \la \nabla d(x_{k}), x - x_k \ra - d(x_{k+1}) - \la \nabla d(x_{k+1}), x - x_{k+1} \ra\\
\\
\refLE{eq-COpt} & \beta_d(x_{k},x)  + d(x_k) + \la \nabla d(x_{k}), x - x_k \ra - d(x_{k+1}) \\
\\
& + \lambda_k \la f'_0(x_k), x - x_{k+1} \ra - \la \nabla d(x_k), x - x_{k+1} \ra\\
\\
= & \beta_d(x_{k},x) + \lambda_k \la f'_0(x_k), x - x_{k+1} \ra - \beta_d(x_k,x_{k+1})\\
\\
= & \beta_d(x_{k},x) - \lambda_k \la f'_0(x_k), x_k - x \ra + \vf_{x_k}(\lambda_k). \QR
\ea
$$

In our method, we will use an {\em indirect control} of the {\em dual} step sizes $\{ \lambda_k \}_{k\geq 0}$, which is based on a predefined sequence of {\em primal} step-size parameters $\{ h_k \}_{k\geq 0}$. As we will prove by Lemma~\ref{lm-Sum}, the convergence of the process can be derived, for example, from the following standard conditions:
\beq\label{eq-Diverge}
\ba{rcl}
\{ h_k \}_{k \geq 0}: \quad \mbox{a) }h_k > 0, \quad \mbox{b) }\sum\limits_{k=0}^{\infty} h_k = \infty.
\ea
\eeq
Then, inequality (\ref{eq-Step}) justifies the choice of the {\em dual} step-size parameter $\lambda_k$ as a solution to the following equation:
\beq\label{eq-Cond}
\mbox{\fbox{$\ba{rcl}
\vf_{x_k}(\lambda_k) & = & \half h_k^2 , \quad k \geq 0
\ea$}}
\eeq

\BE\label{ex-Un}
If $Q = \E$ and $d(x) = \half \| x \|^2_B$, then
$\vf_{\bar x}(\lambda) \refEQ{def-Phi} {\lambda^2 \over 2} (\| f'_0(\bar x) \|^*_B)^2$, and equation (\ref{eq-Cond}) gives us $\lambda_k = h_k /\| f'_0(x_k) \|^*_B$. In this case, method (\ref{met-SGM}) coincides with the classical variant of subgradient method $x_{k+1} = x_k - h_k {f'_0(x_k) \over \| f'_0(x_k) \|^*_B}$. However, for $Q \neq \E$, the rule (\ref{eq-Cond}) allows a proper scaling of the step size by the boundary of feasible set. \qed
\EE

At each iteration of method (\ref{met-SGM}), the corresponding value of $\lambda_k$ can be found by an efficient one-dimensional search procedure, based, for example, on a Newton-type scheme. Since the latter scheme has local quadratic convergence, we make a plausible assumption that it is possible to compute an exact solution of equation (\ref{eq-Cond}). This solution has the following important interpretation (in view of (\ref{eq-DerPhi})):
\beq\label{eq-Inter}
\ba{rcl}
\lambda_k & = & \max\limits_{\lambda \geq 0} \Big\{ \lambda: \; \vf_{x_k}(\lambda) \leq \half h_k^2 \Big\}.
\ea
\eeq
At the same time, we have
$\la \lambda_k f_0'(x_k) + \nabla d(x_{k+1}) - \nabla d(x_k), x_k - x_{k+1} \ra \refGE{eq-COpt} 0$.
Thus
$$
\ba{rcl}
\lambda_k \la f'_0(x_k), x_k - x_{k+1} \ra & \refGE{def-Breg} & \beta_d(x_k,x_{k+1}) + \beta_d(x_{k+1},x_k) \\
\\
& \refGE{eq-BGrow} & \beta_d(x_k,x_{k+1}) +  \half \| x_k - x_{k+1} \|^2.
\ea
$$
Hence
\beq\label{eq-HBound}
\ba{rcl}
\half h_k^2 & \refEQ{eq-Cond} & \vf_{x_k}(\lambda_k) \; \refGE{def-Phi} \; \half \| x_k - x_{k+1} \|^2, \quad k \geq 0.
\ea
\eeq

Our complexity bounds follow from the rate of convergence of the following values:
$$
\ba{rcl}
\delta_k & \Def & \delta_{d_k}(x_k), \quad d_k \; \Def \; {x_k - x_{k+1} \over \| x_k - x_{k+1} \|}.
\ea
$$
\BL\label{lm-Sum}
Let condition (\ref{eq-Cond}) be satisfied. Then for any $N \geq 0$, we have
\beq\label{eq-Rate}
\ba{rcl}
\sum\limits_{k=0}^N h_k \delta_k & \leq & \beta_d(x_0,x^*) + \half \sum\limits_{k=0}^N h_k^2.
\ea
\eeq
\EL
\proof
Indeed, in view of inequality (\ref{eq-Step}), for $r_k = \beta_d(x_k,x^*)$, we have
$$
\ba{rcl}
r_{k+1} & \leq & r_k - \lambda_k \la f_0'(x_k), x_k - x^* \ra + \vf_{x_k}(\lambda_k) \; \refEQ{eq-Cond} \; r_k - \lambda_k \la f_0'(x_k), x_k - x^* \ra + \half h_k^2.
\ea
$$
On the other hand,
$$
\ba{rcl}
\lambda_k \la f_0'(x_k), x_k - x_{k+1} \ra & \refEQ{eq-Cond} & \beta(x_k,x_{k+1}) + \half h_k^2 \; \refGE{eq-BGrow} \; \half \| x_k - x_{k+1} \|^2 + \half h_k^2 \\
\\
& \geq & h_k \| x_k - x_{k+1} \|.
\ea
$$
Hence,
\beq\label{eq-RNext}
\ba{rcl}
r_{k+1} & \leq & r_k - h_k \delta_k + \half h_k^2.
\ea
\eeq
Summing up these inequalities for $k = 0, \dots, N$, we obtain inequality (\ref{eq-Rate}).
\qed

Let us look now at one example of the rate of convergence of method (\ref{met-SGM}) for an objective function from a nonstandard problem class. 
For simplicity, let us measure distances by Euclidean norm $\| x \|_B$ (see  (\ref{def-Euclid})). In this case, we can take $d(x) = \half \| x \|^2_B$ and get $\beta_d(x,y) = \half \| x - y \|^2_B$.
Let us assume that function $f_0(\cdot)$ in problem (\ref{prob-Quasi}) is $p$-times continuously differentiable on $\E$ and its $p$th derivative is Lipschitz-continuous with constant $L_p$. Then we can bound function $\mu(\cdot)$ in (\ref{def-Grow}) as follows:
\beq\label{eq-MuGrow}
\ba{rcl}
\mu(r) & \leq & \sum\limits_{i=1}^p {r^i \over i!} \| D^i f_0(x^*) \|  + {L_p \over (p+1)!} r^{p+1}, \quad r \geq 0,
\ea
\eeq
where all norms for derivatives are induced by $\| \cdot \|_B$.

Let us fix the total number of steps $N \geq 1$ and assume the bound $R_0 \geq \| x_0 - x^* \|_B$ be available. Then, defining the step sizes
\beq\label{eq-SSize}
\ba{rcl}
h_k & = & h \; \Def \; {R_0 \over \sqrt{N+1}}, \quad 0 \leq k \leq N,
\ea
\eeq
we get $\delta^*_N = \min\limits_{0 \leq k \leq N} \delta_k \; \refLE{eq-Rate}{R_0 \over \sqrt{N+1}}$. Hence, in view of inequality (\ref{eq-MuGrow}), we have
\beq\label{eq-FRate}
\ba{rcl}
f_N^* & \Def & \min\limits_{0 \leq k \leq N} f(x_k) \; \leq \; f_0^* + 
\sum\limits_{i=1}^p {1 \over i!} \| D^i f_0(x^*) \| \left( {R_0 \over \sqrt{N+1}} \right)^i  + {L_p \over (p+1)!} \left( {R_0 \over \sqrt{N+1}} \right)^{p+1}.
\ea
\eeq
Note that the first $p$ coefficients in estimate (\ref{eq-Rate}) depend on the local properties of the objective function at the solution, and only the last term employs the global Lipschitz constant for the $p$th derivative. Clearly, we do not need to know the bounds for all these derivatives in order to define the step-size strategy (\ref{eq-SSize}).

\section{Step-size control for max-type convex problems}\label{sc-Conv}
\SetEQ

Let us consider now the following problem of {\em composite optimization}
\beq\label{prob-Comp}
\ba{rcl}
F_* \; \Def \; \min\limits_{x \in \dom \psi} \Big\{ \; F(x) & = & f(x) + \psi(x) \; \Big\},
\ea
\eeq
where $\psi(\cdot)$ is a simple closed convex function, and  
\beq\label{def-FMax}
\ba{rcl}
f(x) & = & \max\limits_{1 \leq i \leq m} \; f_i(x), \quad x \in \dom \psi,
\ea
\eeq
with all $f_i(\cdot)$, $1 \leq i \leq m$, being closed and convex on $\dom \psi$. Denote by
$$
\ba{rcl}
\ell_{\bar x}(x) & = & \max\limits_{1 \leq i \leq m} [ \; f_i(\bar x) + \la f'_i(\bar x), x - \bar x \ra \, ] \; \leq \; f(x), \quad x \in \dom \psi,
\ea
$$
the linearization of function $f(\cdot)$, and by $x_* \in \dom \psi$ an optimal solution of this problem.

Similarly to (\ref{def-Phi}), let us define the following univariate functions:
\beq\label{def-Phi1}
\ba{rcl}
\hat \vf_{\bar x}(\lambda) & = & \max\limits_{x \in \dom \psi} \Big\{ \lambda \Big[ F(\bar x) - \ell_{\bar x}(x) - \psi(x) \Big] - \beta_d(\bar x, x) \Big\}, \quad \lambda \geq 0,
\ea
\eeq
where $\bar x \in \dom \psi$. The unique optimal solution of this optimization problem is denoted by $\hat T_{\bar x}(\lambda)$. Note that function $\hat \vf_{\bar x}(\cdot)$ is convex and continuously differentiable with the following derivative:
\beq\label{eq-DerPhi1}
\ba{rcl}
\hat \vf_{\bar x}'(\lambda) & = & - \Big[ \ell_{\bar x}(\hat T) + \psi(\hat T) - F(\bar x)\Big], \quad \lambda \geq 0,
\ea
\eeq
where $\hat T = \hat T_{\bar x}(\lambda)$.
Thus, $\hat \vf_{\bar x}(0) = 0$, $\hat T_{\bar x}(0) = \bar x$, and $\hat \vf'_{\bar x}(0) = 0$. Hence, $\hat \vf_{\bar x}(\lambda) \geq 0$ for all $\lambda \geq 0$. Let us prove the following variant of Lemma \ref{lm-Step}.
\BL\label{lm-Step1}
Let $\bar x \in \dom \psi$ and $\hat T = \hat T_{\bar x}(\lambda)$ for some $\lambda \geq 0$. Then
\beq\label{eq-Step1}
\ba{rcl}
\beta_d(T,x_*) & \leq & \beta_d(\bar x, x_*) + \lambda \Big[ 
\ell_{\bar x}(x_*) + \psi(x_*) - F(\bar x) \Big]
+ \hat \vf_{\bar x}(\lambda).
\ea
\eeq
\EL
\proof
In view of the first-order optimality condition for the minimization problem in (\ref{def-Phi1}), for all $x \in \dom \psi$, we have
\beq\label{eq-COpt1}
\ba{rcl}
\la \nabla d(\hat T) - \nabla d(\bar x), x - \hat T \ra + \lambda \Big[ \ell_{\bar x}(x)+ \psi(x) \Big] & \geq & \lambda \Big[ \ell_{\bar x}(\hat T)+ \psi(\hat T) \Big].
\ea
\eeq
Therefore, we get
$$
\ba{rcl}
\beta_d(\hat T, x_*) & = & \beta_d(\bar x, x_*) + d(\bar x) + \la \nabla d(\bar x), x_* - \bar x \ra - d(\hat T) - \la \nabla d(\hat T), x^* - \hat T \ra \\
\\
& \refLE{eq-COpt1} & \beta_d(\bar x, x_*) + d(\bar x) + \la \nabla d(\bar x), x_* - \bar x \ra - d(\hat T) \\
\\ 
& & - \la \nabla d(\bar x), x_* - \hat T \ra + \lambda \Big[ \ell_{\bar x}(x_*)+ \psi(x_*) - \ell_{\bar x}(\hat T) -  \psi(\hat T) \Big]\\
\\
& \refEQ{def-Phi1} & \beta_d(\bar x, x_*) + \lambda \Big[ \ell_{\bar x}(x_*)+ \psi(x_*) - F(\bar x) \Big] + \hat \vf_{\bar x}(\hat T). \QR
\ea
$$

In this section, we analyze the following optimization scheme.
\beq\label{met-SGM1}
\ba{|c|}
\hline \\
\quad \mbox{\bf Basic Composite Subgradient Method} \quad\\
\\
\hline \\
\ba{l}
\mbox{{\bf Initialization.}  Choose $x_0 \in \dom \psi$.}\\
\\
\mbox{{\bf $k$th iteration ($k \geq 0$).}}\\
\\
\mbox{a) Choose step-size parameter $\lambda_k > 0$.}\\
\\
\mbox{b) Compute $x_{k+1} = \hat T_{x_k}(\lambda_k)$.}\\
\\
\ea\\
\hline
\ea
\eeq

Suppose that the optimal value $F_*$ is known. In view of convexity of function $f(\cdot)$, we have
$$
\ba{rcl}
\ell_{\bar x}(x_*) + \psi(x_*) - F(\bar x) & \leq & F_* - F(\bar x).
\ea
$$
Hence, for points $\{ x_k \}_{k \geq 0}$, generated by method (\ref{met-SGM1}), we have
\beq\label{eq-Step2}
\ba{rcl}
\beta_d(x_{k+1},x_*) & \refLE{eq-Step1} & \beta_d(x_k,x_*) + \Big(\hat \vf_{x_k}(\lambda_k) - \lambda_k [F(x_k) - F_*] \Big).
\ea
\eeq
This observation explains the following step-size strategy:
\beq\label{eq-Size2}
\ba{c}
\mbox{\fbox{$\lambda_k = \arg\min\limits_{\lambda \geq 0} \Big\{ \; \hat \vf_{x_k}(\lambda) - \lambda [F(x_k) - F_*] \; \Big\}$}}
\ea
\eeq

This strategy has a natural optimization interpretation.
\BL\label{lm-Inter}
Let $\lambda_k$ be defined by (\ref{eq-Size2}). Then
\beq\label{def-THat}
\ba{rcl}
\hat T_{x_k}(\lambda_k) & = & \arg\min\limits_{x \in \dom \psi} \Big\{ \;\beta_d(x_k,x): \; \ell_{x_k}(x) + \psi(x) \leq F_* \; \Big\}.
\ea
\eeq
\EL
\proof
Indeed,
$$
\ba{rl}
& \min\limits_{\lambda \geq 0} \Big\{ \; \hat \vf_{x_k}(\lambda) - \lambda [F(x_k) - F_*] \; \Big\} \\
\\
\refEQ{def-Phi1} & \min\limits_{\lambda \geq 0} \max\limits_{x \in \dom \psi} \Big\{ \; - \lambda \Big[ \ell_{x_k}(x) + \psi(x) -F(x_k) \Big] - \beta_d(x_k,x)  - \lambda [F(x_k) - F_*] \; \Big\}\\
\\
= & \max\limits_{x \in \dom \psi} \min\limits_{\lambda \geq 0} \Big\{ \; - \beta_d(x_k,x) + \lambda \Big[ F_*  - \ell_{x_k}(x) - \psi(x) ]  \Big] \; \Big\}\\
\\
= & - \min\limits_{x \in \dom \psi}\Big\{ \; \beta_d(x_k,x) : \; \ell_{x_k}(x) + \psi(x) \leq F_* \; \Big\}. \QR
\ea
$$

Lemma \ref{lm-Inter} has two important consequences allowing us to estimate the rate of convergence of method (\ref{met-SGM1}) with the step-size rule (\ref{eq-Size2}). Namely, for any $k \geq 0$ we have
\beq\label{eq-DecStep}
\ba{rcl}
\beta_d(x_{k+1},x_*) & \refLE{eq-Step2} & \beta_d(x_k,x_*) - \beta_d(x_k,x_{k+1}),
\ea
\eeq
\beq\label{eq-NewF}
\ba{rcl}
\ell_{x_k}(x_{k+1}) + \psi(x_{k+1}) & \refLE{def-THat} F_*.
\ea
\eeq
Note that method (\ref{met-SGM1}) is not monotone. However, inequality (\ref{eq-DecStep}) gives us a global rate of convergence for the following characteristic:
\beq\label{eq-RateR}
\ba{rcl}
\rho_k & \Def & \min\limits_{0 \leq i \leq k} \| x_i - x_{i+1} \| \; \refLE{eq-BGrow} \;
\min\limits_{0 \leq i \leq k} \sqrt{ 2 \beta_d(x_i,x_{i+1})} \; \refLE{eq-DecStep} \;
{r_0 \over \sqrt{k+1}} , \quad k \geq 0,
\ea
\eeq
where $r_0 = \beta_d(x_0,x_*)$.

\BT\label{th-RateF}
Let all functions $f_i(\cdot)$ have H\"older-continuous gradients:
\beq\label{eq-Hol}
\ba{rcl}
f_i(y) & \leq & f_i(x) + \la f'_i(x), y - x \ra + {L_{\nu} \over 1 + \nu} \| y - x \|^{1+\nu}, \quad x, y \in \dom \psi, \; 1 \leq i \leq p,
\ea
\eeq
where $\nu \in [0,1]$ and $L_{\nu} \geq 0$. Then for the step-size rule (\ref{eq-Size2}) in method (\ref{met-SGM1}) we have
\beq\label{eq-RateF}
\ba{rcl}
F_{k} & \Def & \min\limits_{0 \leq i \leq k} F(x_i) \; \leq \; F_* + {L_{\nu} r_0^{1+\nu} \over 1+\nu}  \left({1 \over k} \right)^{1 + \nu \over 2}, \quad k \geq 1.
\ea
\eeq
\ET
\proof
Indeed, let $\rho_k = \| x_{i_k} - x_{i_k+1} \|$ for some $i_k$, $0 \leq i_k \leq k$. Then for any $k \geq 0$ we have
$$
\ba{rcl}
F_{k+1} \; \leq \; F(x_{i_k + 1}) & \refLE{eq-Hol} & \ell_{x_{i_k}}(x_{i_k + 1}) + {L_{\nu} \over 1 + \nu} \rho_k^{1+\nu} + \psi(x_{i_k+1})\\
\\
& \refLE{eq-NewF} & F_* + {L_{\nu} \over 1 + \nu} \rho_k^{1+\nu}.
\ea
$$
It remains to use inequality (\ref{eq-RateR}).
\qed

Note that the step-size strategy (\ref{eq-Size2}) does not depend on the H\"older parameter $\nu$ in the condition (\ref{eq-Hol}). Hence, the number of iterations of method (\ref{met-SGM1}), (\ref{eq-Size2}), ensuring an $\epsilon$-accuracy in function value, is bounded from above by the following quantity:
\beq\label{eq-UniBound}
\ba{c}
 \inf\limits_{0 \leq \nu \leq 1} \left[ {L_{\nu}  \over (1+\nu) \epsilon} \right]^{2 \over 1 + \nu} \, r_0^2.
\ea
\eeq

To conclude the section, let us show that the step-size rule (\ref{eq-Size2}) can behave much better than the classical one.
\BE\label{ex-OptStep}
Let $\E = R^2$, $\psi(x) = \mbox{\rm Ind} \, Q$ with $Q = \Big\{ x \in \R^2:\; x^{(2)} \leq 0 \Big\}$, $\| x \|^2 = x^T x$,
and $f(x) = \half (x^{(1)})^2 + \half (x^{(2)}-1)^2$. Then $\nabla f(x) = (x^{(1)}, x^{(2)}-1)^T$ and $f_* = 0$. Consider the point $x_k= (x^{(1)}_k, 0)^T$. Then, by the classical rule, we have
$$
\ba{rcl}
x_{k+1} & = & \pi_{Q} \left[ x_k - {f(x_k) - f_* \over \| \nabla f(x_k) \|^2} \nabla f(x_k) \right] \; = \; \pi_{Q} \left[ \left( \ba{c} x^{(1)}_k \\ 0 \ea \right) - {\left(x^{(1)}_k\right)^2 \over 2\left(1+ \left(x^{(1)}_k\right)^2 \right)} \left( \ba{c} x^{(1)}_k \\ - 1 \ea \right) \right].
\ea
$$
Thus, $x^{(2)}_{k+1} = 0$ and $x^{(1)}_{k+1} = x^{(1)}_k - {\left(x^{(1)}_k\right)^3 \over 2\left(1+ \left(x^{(1)}_k\right)^2 \right)}$. This means that $x_k^{(1)} = O(k^{-1/2})$.

On the other hand, the rule (\ref{def-THat}) as applied to the same point $x_k$ defines $x = x_{k+1}$ as an intersection of two lines
$$
\ba{rcl}
\Big\{x \in \R^2: \; \la \nabla f(x_k), x_k - x \ra = f(x_k) - f_* \Big\}, \quad \Big\{x\in \R^2: \; x^{(2)} = 0 \Big\}.
\ea
$$
This means that $x^{(1)}_k \left(x^{(1)}_k - x^{(1)}_{k+1}\right) = \half \left(x^{(1)}_k\right)^2$, and we get the linear rate of convergence to the optimal point $x_* = 0$. 
\qed
\EE

\section{Step-size control for problems with smooth \\ strongly convex components}\label{sc-RelSmooth}
\SetEQ

A linear rate of convergence, demonstrated by method (\ref{met-SGM1}) in the Example \ref{ex-OptStep}), provides us with motivation to look at behavior of this method on smooth and strongly convex problems. Let us introduce a Euclidean metric by (\ref{def-Euclid}) and define $d(x) = \half \| x \|^2_B$, $x \in \E$. In this case,
$$
\ba{rcl}
\beta_d(x,y) & = & \half \| x - y \|^2_B, \quad x, y \in \E.
\ea
$$

Suppose that all functions $f_i(\cdot)$, $i = 1, \dots, m$ in (\ref{def-FMax}) are continuously differentiable and have Lipschitz-continuous gradients with the same constant $L_f \geq 0$:
\beq\label{def-RSmooth}
\ba{rcl}
f_i(y) & \leq & f_i(x) + \la \nabla f_i(x), y - x \ra + \half L_f \| x - y \|^2_B, \quad x, y \in \dom \psi .
\ea
\eeq
In the notation of Section \ref{sc-Conv}, these inequalities imply the following upper bound for the objective function of problem (\ref{prob-Comp}):
\beq\label{eq-FCUp}
\ba{rcl}
F(y) & \leq & \ell_x(y) + \psi(y) + \half L_f \| x - y \|^2_B, \quad x, y \in \dom \psi.
\ea
\eeq
Hence, for the sequence of points $\{ x_k \}_{k \geq 0}$, generated by method  (\ref{met-SGM1}), we have
\beq\label{eq-NFBound}
\ba{rcl}
F(x_{k+1}) & \leq & \ell_{x_k}(x_{k+1}) + \psi(x_{k+1}) + \half L_f \| x_k - x_{k+1} \|^2_B \\
\\
& \refLE{eq-NewF} & F_* + \half L_f \| x_k - x_{k+1} \|^2_B.
\ea
\eeq
Thus, we can prove the following statement.
\BT\label{th-Smooth}
Under condition (\ref{def-RSmooth}), the rate of convergence of method (\ref{met-SGM1}) can be estimated as follows:
\beq\label{eq-CRate1}
\ba{rcl}
{1 \over T} \sum\limits_{k=1}^{T} F(x_k) - F_* & \leq & {1 \over 2T} L_f \| x_0 - x_* \|^2_B, \quad T \geq 1.
\ea
\eeq
If in addition functions $f_i(\cdot)$, $i = 1, \dots, m$, are  strongly convex:
\beq\label{def-SConv}
\ba{rcl}
f_i(y) & \geq & f_i(x) + \la \nabla f_i(x), y - x \ra + \half \mu_f \| x - y \|^2_B, \quad x, y \in \dom \psi, 
\ea
\eeq
with $\mu_f > 0$, then the rate of convergence is linear:
\beq\label{eq-SCRate1}
\ba{rcl}
\| x_k - x_* \|^2_B & \leq & \left({L_f \over \mu_f + L_f}\right)^k \| x_0 - x_* \|^2_B, \quad k \geq 0.
\ea
\eeq
\ET
\proof
Indeed, substituting inequality (\ref{eq-NFBound}) into relation (\ref{eq-DecStep}), we have
\beq\label{eq-DecSF}
\ba{rcl}
\half \| x_{k+1} - x_* \|^2_B & \leq & \half \| x_k - x_* \|^2_B - {1 \over L_f} [ F(x_{k+1}) - F_*], \quad k \geq 0.
\ea
\eeq
Summing up these inequalities for $k = 0, \dots, T-1$, we get inequality (\ref{eq-CRate1}).

If in adition, the conditions (\ref{def-SConv}) are satisfied, then
\beq\label{eq-FCDown}
\ba{rcl}
F(y) & \geq & \ell_x(y) + \psi(y) + \half \mu_f \| x - y \|^2_B, \quad x, y \in \dom \psi.
\ea
\eeq
Note that the first-order optimality conditions for problem (\ref{prob-Comp}) can be written in the following form:
$$
\ba{rcl}
\ell_{x_*}(y) + \psi(y) & \geq & F_*, \quad \forall y \in \dom \psi.
\ea
$$
Therefore, inequality (\ref{eq-FCDown}) implies that
$$
\ba{rcl}
F(y) - F_* & \geq & \half \mu_f \| y - x_* \|^2_B, \quad \forall y \in \dom \psi.
\ea
$$
Thus, $\| x_{k+1} - x_* \|^2_B \refLE{eq-DecSF} {L_f \over \mu_f + L_f} \| x_k - x_* \|^2_B$, and we get inequality (\ref{eq-SCRate1}).
\qed

\section{Convex minimization with max-type composite constraint}\label{sc-Cons}
\SetEQ

Let us show that both step-size strategies described in Sections \ref{sc-Quasi} and \ref{sc-Conv} can be unified in one scheme for solving constrained optimization problems. In this section, we deal with the problem in the following {\em semi-composite} form:
\beq\label{prob-Semi}
\ba{rcl}
\min\limits_{x \in Q} \Big\{ f_0(x): \; F(x) \Def f(x) + \psi(x) \leq 0 \Big\}
\ea
\eeq
where $\psi(\cdot)$ is a simple closed convex function, $Q \subseteq \dom \psi$ is a closed convex set, and
$$
\ba{rcl}
f(x) & = & \max\limits_{1 \leq i \leq m} f_i(x),
\ea
$$
with all functions $f_i(\cdot)$, $i = 0, \dots, m$, being closed and convex on $\dom \psi$. We assume $Q$ to be bounded: 
\beq\label{eq-DBound}
\ba{rcl}
\beta_d(x,y) & < & D, \quad \forall x, y \in Q.
\ea
\eeq

In order to solve problem (\ref{prob-Semi}), we propose a method, which combines two different types of iterations. One of them improves the {\em feasibility} of the current point, and the second one improves its {\em optimality}.

Iteration of the first type is based on the machinery developed in Section \ref{sc-Conv} with the particular value $F_* = 0$. It is applied to some point $x_k \in Q$.
\beq\label{eq-LInfeas}
\mbox{\fbox{$\ba{l}\\
\quad \mbox{Compute $\lambda_k = \arg\min\limits_{\lambda \geq 0} \Big\{ \hat \vf_{x_k}(\lambda) - \lambda F(x_k) \Big\}$. Set $x_{k+1} \refEQ{def-THat} \hat T_{x_k}(\lambda_k)$.} \quad\\
\\
\ea$}}
\eeq
Note that for $F(x_k) \leq 0$, we have $\lambda_k = 0$ and $\hat T_{x_k}(\lambda_k) = x_k$.

For iteration $k$ of the second type, we need to choose a primal step size bound $h_k > 0$. Then, at the test point $y_k \in Q$, we define the function $\vf_{y_k}(\cdot)$ by (\ref{def-Phi}) and apply the following rule:
\beq\label{def-HConst}
\mbox{\fbox{$\ba{l}
\\
\quad\mbox{Find $\lambda_k$ from equation
$\vf_{y_k}(\lambda_k) = \half h_k^2$. Set $x_{k+1} = T_{y_k}(\lambda_k)$.}\quad \\
\\
\ea$}}
\eeq

Since in both rules parameters $\lambda_k$ are functions of the test points, we will use shorter notations $T(y_k)$ and $\hat T(x_k)$. Consider the following optimization scheme.

\beq\label{met-SGM2}
\ba{|c|}
\hline \\
\quad \mbox{\bf Double-Step Subgradient Method for Semi-Composite Problem} \quad\\
\\
\hline \\
\ba{l}
\mbox{{\bf Initialization.}  Choose $x_0 \in \dom \psi$ and sequence of steps ${\cal H} = \{ h_k \}_{k \geq 0}$.}\\
\\
\mbox{{\bf $k$th iteration ($k \geq 0$).} {\bf 1.} Compute $y_k=\hat T(x_k)$.}\\
\\
\mbox{{\bf 2.} If $\beta_d(x_k,y_k) \geq \half h_k^2$, then a) $x_{k+1} = y_k$. Else, b) compute $x_{k+1} = T(y_k)$.}\\
\\
\ea\\
\hline
\ea
\eeq

Thus, method (\ref{met-SGM2}) is defined by a sequence of primal step bounds ${\cal H} = \{ h_k \}_{k \geq 0}$. However, since in (\ref{met-SGM2}) we apply a switching strategy, it is impossible to say in advance what will be the type of a particular $k$th iteration. Therefore, as compared with the classical conditions (\ref{eq-Diverge}), we need additional regularity assumptions on ${\cal H}$.

It will be convenient to relate this sequence with another sequence of {\em scaling coefficients} ${\cal T} = \{ \tau_k \}_{k\geq 0}$, satisfying the {\em second-order divergence condition} (compare with (\ref{eq-Diverge})).
\beq\label{eq-TGrow}
\ba{|c|}
\hline \\
\mbox{\bf Second-Order Divergence Condition}\\
\\
\hline \\
\ba{rl}
\mbox{\bf a)} & \mbox{$\tau_k \geq \tau_{k+1} > 0$ for any $k \geq 0$.} \quad \mbox{ {\bf b)} $\sum\limits_{k=0}^{\infty} \tau_k^2 = + \infty$.}
\ea\\
\\
\hline
\ea
\eeq

Note that condition (\ref{eq-TGrow}) ensures $\sum\limits_{k=0}^{\infty} \tau_k = + \infty$. Thus, it is stronger than (\ref{eq-Diverge}). In order to transform the sequence ${\cal T}$ into {\em convergence rate} of some optimization process, we need to introduce the following characteristic.
\BD\label{def-DDelay}
For a sequence ${\cal T}$, the integer-valued function $a(k)$, $k \geq 0$, is called the {\em divergence delay} (of degree two) if $a(k) \geq 0$ is the minimal integer value such that
\beq\label{def-DD}
\ba{rcl}
\sum\limits_{i=k}^{k+a(k)} \tau_i^2 & \geq & 1, \quad k \geq 0.
\ea
\eeq
\ED
Clearly, condition (\ref{eq-TGrow})$_b$ ensures that all values $a(k)$, $k \geq 0$, are well defined. At the same time, from (\ref{eq-TGrow})$_a$, we have $a(k+1) \geq a(k)$ for any $k \geq 0$.

Let us give two important examples of such sequences.
\BE\label{eq-T2}
{\bf a)} Let us fix an integer $N \geq 0$ and define
\beq\label{def-T0}
\ba{rcl}
\tau_k & = & {1 \over \sqrt{N+1}}, \quad k \geq 0.
\ea
\eeq
Then condition (\ref{eq-TGrow})$_a$ is valid and 
$\sum\limits_{i=k}^{k+a(k)} \tau_i^2 = {a(k)+1 \over N+1}$. Thus, $a(k)= N$.

{\bf b)} Consider the following sequence:
\beq\label{def-T1}
\ba{rcl}
\tau_k & = & \sqrt{2 \over k+1}, \quad k \geq 0.
\ea
\eeq
For $k \geq 1$, denote by $S_k = \sum\limits_{i=0}^{k-1} \tau_i^2$.
Then $S_1 = 2$, $S_2 = 3$, and the difference 
$$
\ba{rcl}
S_{2k} - S_k & = & \sum\limits_{i=k}^{2k-1} \tau_i^2
\ea
$$
is monotonically increasing in $k \geq 1$. Thus, $a(k) \leq k-1$ for all $k \geq 1$.
\qed	
\EE

Let us analyze performance of method (\ref{met-SGM2}) with an appropriately chosen sequence ${\cal H}$. Namely, let us  choose
\beq\label{def-HStep2}
\ba{rcl}
h_k & = &  \sqrt{2D} \; \tau_k, \quad k \geq 0,
\ea
\eeq
where the sequence ${\cal T}$ satisfies condition (\ref{eq-TGrow}) and $D$ is taken from (\ref{eq-DBound}). 

We are interested only in the values of  objective function $f_0(\cdot)$ computed at the points with small values of the functional constraint $F(\cdot)$. As we will see, these points are involved in Step~2b). For the total number of steps $N \geq 1 + a(0)$, denote
\beq\label{def-KN}
\ba{rcl}
k(N) & = & \max\Big\{ k \geq 0: \; k + a(k) \leq N-1\Big\},\\
\\
{\cal F}_N & = & \Big\{ k: \; k(N) \leq k \leq N-1, \; y_k \stackrel{2b)}{\longrightarrow} x_{k+1} \Big\}.
\ea
\eeq

Let us define the following directional proximity measures:
$$
\ba{rclll}
\delta_k & = & \delta_{d_k}(y_k) & \mbox{with $\; d_k = {y_k - x_{k+1} \over \| y_k - x_{k+1} \|}$}, & k \in {\cal F}_N.
\ea
$$
We are interested in the rate of convergence to zero of the following characteristic:
$$
\ba{rcl}
\delta_N^* & = & \min\limits_{k \in {\cal F}_N} \delta_k.
\ea
$$

\BT\label{th-Main}
For any $N \geq 1 + a(0)$, the number $k(N)\geq 0$ is well defined and $\delta_N^* < h_{k(N)}$.

\noindent
Moreover, if all functions $f_i(\cdot)$, $i = 1, \dots, m$, have H\"older-continuous gradients on $Q$ with parameter $\nu \in [0,1]$ and constant $L_{\nu}>0$, then, for any $k \in {\cal F}_N$, we have
\beq\label{eq-FSmall}
\ba{rcl}
F(y_k) & \leq & {L_{\nu} \over 1 + \nu} h_{k(N)}^{1+\nu}.
\ea
\eeq
\ET
\proof
Let us bound the distances $r_k = \beta_d(x_k,x_*)$ for $k \geq k(N)$. If $k \not\in {\cal F}_N$, then
$$
\ba{rcl}
r_{k+1} & = & \beta_d(y_k ,x_* ) \; \refLE{eq-DecStep} \; r_k - \beta(x_k,y_k ) \; \refLE{met-SGM2} \; r_k - \half h_k^2.
\ea
$$
If $k \in {\cal F}_N$, then
$$
\ba{rcl}
r_{k+1} & \refLE{eq-RNext} & \beta(y_k,x_*) - h_k \delta_k + \half h_k^2 \; \refLE{eq-DecStep} \; r_k - h_k \delta_k + \half h_k^2 \; \leq \;r_k - h_k \delta^*_N + \half h_k^2 .
\ea
$$
Summing up these inequalities for $k = k(N), \dots N-1$, we obtain
$$
\ba{rcl}
r_{N} & < & D - \half \sum\limits_{k \not\in {\cal F}_N} h_k^2 
+ \half \sum\limits_{k \in {\cal F}_N} ( h_k^2 - 2 h_k \delta^*_N)\\
\\
& = & D - \half \sum\limits_{k=k(N)}^{N-1} h_k^2 
+ \sum\limits_{k \in {\cal F}_N} ( h_k^2 - h_k \delta^*_N)\\
\\
& \stackrel{(\ref{def-DD}),(\ref{def-HStep2})}{<} & \sum\limits_{k \in {\cal F}_N} h_k ( h_k - \delta^*_N) \; \stackrel{(\ref{eq-TGrow})_a}{\leq} \; ( h_{k(N)} - \delta^*_N) \sum\limits_{k \in {\cal F}_N} h_k .
\ea
$$
Thus, $\delta^*_N < h_{k(N)}$.
Finally, inclusion $k \in {\cal F}_N$ implies 
$$
\ba{rcl}
\half \| y_k - x_k \|^2 & \refLE{eq-BGrow} & \beta_d(x_k,y_k) \; \leq \; \half h_k^2 \; \stackrel{(\ref{eq-TGrow})_a}{\leq} \; \half h_{k(N)}^2. 
\ea
$$
Since $F_* = 0$, we get
$$
\ba{rcl}
F(y_k) & \leq & \ell_{x_k}(y_k) + {L_\nu \over 1+\nu} \| y_k - x_k \|^{1+\nu} + \psi(y_k)\\
\\
& \refLE{eq-NewF} & {L_\nu \over 1+\nu} \| y_k - x_k \|^{1+\nu} \; \leq \; {L_\nu \over 1+\nu} h_{k(N)}^{1+\nu},
\ea
$$
and this is inequality (\ref{eq-FSmall}).
\qed

As a straightforward consequence of Theorem \ref{th-Main}, we have the following rate of convergence in function value:
\beq\label{eq-Rate2}
\ba{rcl}
\min\limits_{x \in {\cal F}_N} f_0(y_k) & \leq & f^*_0 + \mu(h_{k(N)}),
\ea
\eeq
where the function  $\mu(\cdot)$ is defined by (\ref{def-Grow}). 

Thus, the actual rate of convergence of method (\ref{met-SGM2}) depends on the rate of convergence of sequence ${\cal T}$ and the magnitude of divergence delay. For example, for the choice (\ref{def-T1}), in view of inequality $a(k) \leq k-1$, we have
$$
\ba{rcl}
k(N) & = & \max\Big\{ k : \; k + a(k) \leq N-1\Big\} \; \geq \; 
\max\Big\{ k : \;2k-1 \leq N-1\Big\}.
\ea
$$
Hence, for $N = 2M$, the choice (\ref{def-T1}) ensures $k(N) \geq M = N/2$, and we have 
\beq\label{eq-HRate}
\ba{rcl}
h_{k(N)} & \leq & h_M \; = \; {2 D^{1/2} \over \sqrt{M + 1}}. 
\ea
\eeq
It is interesting that the rate of convergence (\ref{eq-FSmall}) for the constraints can be higher than the rate (\ref{eq-Rate2}) for the objective function.

\section{Approximating Lagrange Multipliers, I}\label{sc-Lagrange1}
\SetEQ

Despite to the good convergence rate (\ref{eq-Rate2}), when the number of functional components $m$ in problem (\ref{prob-Semi}) is big, the implementation of one iteration (\ref{eq-LInfeas}) can be very expensive. In this section, we consider a simpler switching strategy for solving convex optimization problems with potentially many functional constraints. Our method is also able to approximate the corresponding optimal Lagrange multipliers.

Consider the following {\em constrained} optimization problem:
\beq\label{prob-FConst}
\ba{rcl}
f_0^* \; = \; \min\limits_{x \in Q} \Big\{ f_0(x) : \; f_i(x)  \leq 0, \; i = 1, \dots, m \Big\},
\ea
\eeq
where all functions $f_i(\cdot)$ are closed and convex, $i = 0, \dots, m$, and set $Q$ is closed, convex, and bounded. Denote by $x^*$ one of its optimal solutions. Sometimes we use notation $\bar f(x) = (f_1(x), \dots, f_m(x))^T \in \R^m$. Denote ${\cal F} = \{ x \in Q: \; \bar f(x) \leq 0 \}$.

We assume that functions $f_i(\cdot)$ are subdifferentiable on $Q$ and it is possible to compute their subgradients with uniformly bounded norms:
\beq\label{eq-GBound}
\ba{rcl}
\| f'_i(x) \|_* & \leq & M_i, \quad x \in Q, \quad i = 0, \dots, p.
\ea
\eeq
Denote $\bar M = (M_1, \dots, M_m)^T \in \R^m$.

For the set $Q$, we assume existence of a prox-function $d(\cdot)$ defining the corresponding Bregman distance $\beta_d(\cdot,\cdot)$. In our methods, we need to know a constant $D$ such that
\beq\label{eq-DHBound}
\ba{rcl}
\beta_d(x,y) & < & D, \quad \forall x, y \in Q.
\ea
\eeq

Let us introduce the Lagrangian 
$$
\ba{rcl}
{\cal L}(x,\bar \lambda) & = & f_0(x) + \la \bar \lambda, \bar f(x) \ra, \quad \bar \lambda = (\lambda^{(1)}, \dots, \lambda^{(m)})^T \in \R^m_+, 
\ea
$$
and the dual function $\phi(\bar \lambda) = \min\limits_{x \in Q} {\cal L}(x,\bar \lambda)$. By Sion's theorem, we know that
\beq\label{prob-LDual}
\ba{rcl}
\sup\limits_{\bar \lambda \in \R^m_+} \phi(\bar \lambda) & = & f_0^*.
\ea
\eeq

Our first method is based on two operations presented in Sections \ref{sc-Quasi} and \ref{sc-Conv}. For defining iterations of the first type, we need functions $\vf_{i,x}(\lambda)$ with $\lambda \geq 0$, parameterized by $x \in Q$, and defined as follows:
\beq\label{def-NPhi}
\ba{rcl}
\vf_{i,x}(\lambda) & = & \max\limits_{T \in Q} \Big\{ \lambda \la f'_i(x), x - T \ra - \beta_d(x,T) \Big\}, \quad i = 0, \dots , m.
\ea
\eeq
An appropriate value of $\lambda$ can be found from the equation
\beq\label{eq-HLambda}
\ba{rcl}
\vf_{i,x_k}(\lambda) & = & \half h_k^2,
\ea
\eeq
and used for setting $x_{k+1} = T_{i,x_k}(\lambda)$, the optimal solutions of  problem (\ref{def-NPhi}). In this section, we perform this iteration only for the objective function ($i=0$). The possibility of using the steps $T_{i,x_k}(\lambda)$ for inequality constraints is analyzed in Sections \ref{sc-Lagrange2} and \ref{sc-NoBound}.

The iteration of the second type is trying to improve feasibility of the current point $x_k \in Q$ (compare with (\ref{eq-LInfeas})). It needs computation of all Bregman projections
\beq\label{def-TLambda}
\ba{rcl}
\hat T_i(x_k) & \Def & \arg\min\limits_{T \in Q} \Big\{ \beta_d(x_k,T): \;  f_i(x_k) + \la f'_i(x_k), T - x_k \ra \leq 0 \Big\}, 
\ea
\eeq
for indexes $i = 1, \dots, m$. The first-order optimality condition for  problem (\ref{def-TLambda}) is as follows:
\beq\label{eq-TLambda}
\ba{rcl}
\la\nabla d(T) - \nabla d(x_k) + \lambda f'_i(x_k), x - T \ra  & \geq & 0, \quad \forall x \in Q,
\ea
\eeq
where $T = \hat T_i(x_k)$ and $\lambda = \lambda_i(x_k) \geq 0$ is the optimal Lagrange multiplier for the linear inequality constraint in (\ref{def-TLambda}). We assume that this multiplier can be also computed. Note that for $x = x_k$ we have
$$
\ba{rcl}
\lambda_i(x_k) \la f'_i(x_k), x_k - T \ra  & \geq & \beta_d(x_k,T) + \beta_d(T,x_k) \; \refGE{eq-BGrow} \; \| x_k - T \|^2.
\ea
$$
Hence, if $T \neq x_k$, then $\lambda_i(x_k) > 0$.

Consider the following optimization scheme.
\beq\label{met-SGM3}
\ba{|c|}
\hline \\
\quad \mbox{\bf Subgradient Method for Problem (\ref{prob-FConst})} \quad\\
\\
\hline \\
\ba{l}
\mbox{{\bf Initialization.}  Choose $x_0 \in Q$ and a sequence of step bounds ${\cal H} = \{ h_k \}_{k \geq 0}$.}\\
\\
\mbox{{\bf $k$th iteration ($k \geq 0$).} {\bf 1.} Compute $r_k^i= \beta_d(x_k,\hat T_i(x_k))$, $1 \leq i \leq m$.}\\
\\
\mbox{{\bf 2.} If $\exists \; i \in \{1 \dots m\}: r_k^{i} \geq \half h^2_k$, then set $i_k = i$, $\lambda_k = \lambda_i(x_k)$, $x_{k+1} = \hat T_{i}(x_k)$.}\\
\\
\mbox{{\bf 3.} Else, set $i_k = 0$, compute $\lambda_k$ by (\ref{eq-HLambda}) and set $x_{k+1} = T_{x_k}(\lambda_k)$.}\\
\\
\ea\\
\hline
\ea
\eeq

Let us prove that method (\ref{met-SGM3}) can find an approximate solution of the primal-dual problem (\ref{prob-FConst}), (\ref{prob-LDual}). Let us choose a sequence ${\cal T}$ satisfying condition (\ref{eq-TGrow}) and define 
\beq\label{def-HStep3}
\ba{rcl}
h_k & = & \sqrt{2D} \; \tau_k, \quad k \geq 0,
\ea
\eeq
where $D$ satisfies inequality (\ref{eq-DHBound}). Then, for the number of steps $N \geq 1 + a(0)$, let us define function $k(N)$ by the first equation in (\ref{def-KN}). Now we can define the following objects:
\beq\label{def-Sigma}
\ba{rcl}
{\cal A}_i(N) & = & \{k: \; k(N) \leq k \leq N-1, \; i_k = i \}, \\
\\
\sigma_i(N) & = & \sum\limits_{k \in {\cal A}_i(N)} \lambda_k, \quad i = 0, \dots, m,\\
\\
\lambda^*_i(N) & = & \sigma_i(N)/\sigma_0(N), \quad i = 1, \dots, m.
\ea
\eeq
Clearly, the vectors of dual multipliers $\bar \lambda_*(N) = (\lambda^*_1(N), \dots, \lambda^*_m(N))$ are defined only if ${\cal A}_0(N) \neq \emptyset$. As usual, the sum over an empty set of iterations is assumed to be zero.

\BT\label{th-Main2}
Let the sequence of step bounds ${\cal H}$ in method (\ref{met-SGM3}) be defined by (\ref{def-HStep3}).
Then for any $N \geq 1 + a(0)$, we have ${\cal A}_0(N) \neq \emptyset$.
Moreover, if all functions $f_i(\cdot)$, $0 = 1, \dots, m$, satisfy (\ref{eq-GBound}), then
\beq\label{eq-GSmall}
\ba{rcl}
\max\limits_{1 \leq i \leq m} {1 \over M_i} f_i(x_k) & \leq & h_{k(N)}, \quad k \in {\cal A}_0(N),
\ea
\eeq
\beq\label{eq-DGSmall}
\ba{rcl}
{1 \over \sigma_0(N)} \sum\limits_{k \in {\cal A}_0(N)} \lambda_k f_0(x_k)  & \leq & \phi(\bar \lambda^*(N)) + M_0 h_{k(N)}.
\ea
\eeq
\ET
\proof
Indeed, we have
$$
\ba{rcl}
B_N & \Def & \max\limits_{x \in Q} \Big[  \sum\limits_{k \in {\cal A}_0(N)} \lambda_k f_0(x_k)  - \sigma_0(N) f_0(x) - \sum\limits_{i=1}^m \sigma_i(N) f_i(x) \Big]\\
\\
& = & \max\limits_{x \in Q} \Big[  \sum\limits_{k \in {\cal A}_0(N)} \lambda_k [f_0(x_k)  - f_0(x)] - \sum\limits_{i=1}^m \sigma_i(N) f_i(x) \Big].
\ea
$$
Let us fix some $x \in Q$ and denote $r_k(x) = \beta_d(x_k,x)$. Then, for $k \in {\cal A}_0(N)$, we have
\beq\label{eq-F0}
\ba{rcl}
\lambda_k [f_0(x_k)  - f_0(x)] & \leq &  \lambda_k \la f'_0(x_k), x_k - x \ra.
\ea
\eeq
At the same time, from the first-order optimality condition for problem (\ref{def-NPhi}), we have
$$
\ba{rl}
0 \; \leq & \la \nabla d(x_{k+1}) - \nabla d(x_k), x - x_{k+1} \ra + \lambda_k \la f'_0(x_k), x - x_{k+1} \ra \\
\\
\refLE{eq-F0} & \la \nabla d(x_{k+1}) - \nabla d(x_k), x - x_{k+1} \ra - \lambda_k [f_0(x_k) - f_0(x)] + \lambda_k \la f'_0(x_k), x_k - x_{k+1} \ra  \\
\\
= &  \la \nabla d(x_{k+1}) - \nabla d(x_k), x - x_{k+1} \ra - \lambda_k [f_0(x_k) - f_0(x)] + \vf_{0,x_k}(\lambda_k) + \beta_d(x_k,x_{k+1}).
\ea
$$
Hence,
$$
\ba{rl}
& r_{k+1}(x) - r_k(x) \; = \; \la \nabla d(x_k) - \nabla d(x_{k+1}), x - x_{k+1} \ra - \beta_d(x_k,x_{k+1})\\
\\
\leq & - \lambda_k [f_0(x_k) - f_0(x)] + \vf_{0,x_k}(\lambda_k) \; \refEQ{eq-HLambda} \; - \lambda_k [f_0(x_k) - f_0(x)] + \half h_k^2.
\ea
$$

On the other hand,
$$
\ba{rcl}
\sum\limits_{i=1}^m \sigma_i(N) f_i(x) & \refEQ{def-Sigma} & \sum\limits_{i=1}^m \sum\limits_{k \in {\cal A}_i(N)} \lambda_k f_i(x) \; = \; 
\sum\limits_{k \not\in {\cal A}_0(N)} \lambda_{i_k}(x_k) f_{i_k}(x)\\
\\
& \geq & \sum\limits_{k \not\in {\cal A}_0(N)} \lambda_{i_k}(x_k) [f_{i_k}(x_k) + \la f'_{i_k}(x_k), x - x_k \ra ]\\
\\
& \refEQ{def-TLambda} & \sum\limits_{k \not\in {\cal A}_0(N)} \lambda_{i_k}(x_k) [\la f'_{i_k}(x_k), x - x_{k+1} \ra ] \\
\\
& \refGE{eq-TLambda} & \sum\limits_{k \not\in {\cal A}_0(N)} \la \nabla d(x_{k+1}) - \nabla d(x_k), x_{k+1} - x \ra\\
\\
& = & \sum\limits_{k \not\in {\cal A}_0(N)} \Big[ r_{k+1}(x) - r_{k}(x) + \beta_d(x_k,x_{k+1}) \Big] \\
\\
& \stackrel{(\ref{met-SGM3})_2}{\geq} & \sum\limits_{k \not\in {\cal A}_0(N)} \Big[ r_{k+1}(x) - r_{k}(x) + \half h_k^2 \Big].
\ea
$$
Therefore,
$$
\ba{rcl}
B_N & \leq & \max\limits_{x \in Q} \Big[ \sum\limits_{k \in {\cal A}_0(N)} [r_k(x) - r_{k+1}(x) + \half h_k^2] - \sum\limits_{k \not\in {\cal A}_0(N)} [ r_{k+1}(x) - r_{k}(x) + \half h_k^2 ] \; \Big]\\
\\
& = & \max\limits_{x \in Q} \Big[ r_{k(N)}(x) - r_N(x) + \half \sum\limits_{k \in {\cal A}_0(N)} h_k^2 - \half \sum\limits_{k \not\in {\cal A}_0(N)} h_k^2 \Big] \\
\\
& < & D - \half \sum\limits_{k=k(N)}^{N-1} h_k^2 + \sum\limits_{k \in {\cal A}_0(N)} h_k^2 \; \refLE{def-HStep3} \; \sum\limits_{k \in {\cal A}_0(N)} h_k^2 \; \stackrel{(\ref{eq-TGrow})_a}{\leq} \;  h_{k(N)} \sum\limits_{k \in {\cal A}_0(N)} h_k. 
\ea
$$

Let us assume now that ${\cal A}_0(N) = \emptyset$. Then $B_N < 0$, and this is impossible since the point $x_*$ is feasible. Thus, we have proved that ${\cal A}_0(N) \neq \emptyset$.

Further, for any $k \in {\cal A}_0(N)$ we have
$$
\ba{rcl}
\lambda_k \la f'_0(x_k), x_k - x_{k+1} \ra  & \refEQ{eq-HLambda} & \beta_d(x_k,x_{k+1}) + \half h_k^2.
\ea
$$
Hence, $r_k^0 = \| x_{k+1} - x_k \| > 0$ and we conclude that
$$
\ba{rcl}
\lambda_k M_0 & \refGE{eq-GBound} \lambda_k \| f'_0(x_k) \|_* \; \geq {\lambda_k \over r_k^0} \la f'_0(x_k), x_k - x_{k+1} \ra  \; \refGE{eq-BGrow} \; {1 \over 2 r_k^0} \Big[ (r_k^0)^2 + h_k^2 \Big] \; \geq \; h_k.
\ea
$$
Thus, we have proved that $B_N \leq M_0 \sigma_0(N) h_{k(N)}$. Dividing this inequality by $\sigma_0(N) > 0$, we get the bound (\ref{eq-DGSmall}).

Finally, note that for any $k \in {\cal A}_0(N)$ we have $\| \hat T_i(x_k) - x_k \| \leq h_k$ for all $i$, $1 \leq i \leq m$. Note that either $f_i(x_k) \leq 0$, or $f_i(x_k) > 0$ and $f_i(x_k) + \la f'_i(x_k), \hat T_i(x_k) - x_k) \ra = 0$. In the latter case, we have
$$
\ba{rcl}
f_i(x_k) & = & \la f'_i(x_k), x_k - \hat T_i(x_k)\ra \; \refLE{eq-GBound} \; M_i h_k \; \leq \; M_i h_{k(N)}.
\ea
$$
Thus, inequality (\ref{eq-GSmall}) is proved.
\qed

Note that for the choice of scaling sequence (\ref{def-T1}), method (\ref{met-SGM2}) is globally convergent. In this case, the computation of Lagrange multipliers by (\ref{def-Sigma}) for all values of $N$, requires storage of all coefficients $\{ \lambda_k \}_{k \geq 0}$. This inconvenience can be avoided if we decide to accumulate the sums for Lagrange multipliers only starting from the moments $k(N_q)$ with $N_q = 2^q$, $q \geq 1$. Then the method (\ref{met-SGM2}) will be allowed to stop only at the moments $2N_q$.

\section{Approximating Lagrange Multipliers, II}\label{sc-Lagrange2}
\SetEQ

Method (\ref{met-SGM3}) has one hidden drawback. If $i_k \geq 1$, then
$$
\ba{rcl} 
0 & \refEQ{def-TLambda} & f_{i_k}(x_k) + \la f'_{i_k}(x_k), x_{k+1} - x_k \ra \; \leq \; f_{i_k}(x_{k+1}).
\ea
$$
Thus, this scheme most probably generates infeasible approximations of the optimal point, which violate some of the functional constraints. In order to avoid this tendency, we propose a scheme which uses for both types of iterations (improving either feasibility or optimality) the same step-size rule (\ref{eq-HLambda}).
\beq\label{met-SGM4}
\ba{|c|}
\hline \\
\quad \mbox{\bf Fixed-Step Switching Subgradient Method for Problem (\ref{prob-FConst})} \quad\\
\\
\hline \\
\ba{l}
\mbox{{\bf Initialization.}  Choose $x_0 \in Q$ and a sequence of step bounds ${\cal H} = \{ h_k \}_{k \geq 0}$.}\\
\\
\mbox{{\bf $k$th iteration ($k \geq 0$).} {\bf 1.} Compute $T_{i,x_k}( \lambda_{i,k})$ by (\ref{eq-HLambda}), $1 \leq i \leq m$.}\\
\\
\mbox{{\bf 2.} If $\exists \; i \in \{1 \dots m\}: \lambda_{i,k} f_i(x_k) \geq  h^2_k$, then set $i_k = i$ and $x_{k+1} = T_{i,x_k}(\lambda_{i,k})$.}\\
\\
\mbox{{\bf 3.} Else, set $i_k = 0$, compute $T_{0,x_k}( \lambda_{0,k})$ by (\ref{eq-HLambda}), and set $x_{k+1} = T_{0,k}(\lambda_{0,k})$.}\\
\\
\mbox{{\bf 4.} Set $\lambda_k = \lambda_{i_k,k}$.}\\
\\
\ea\\
\hline
\ea
\eeq

Thus, for both types of iterations, we use the same step-size strategy (\ref{eq-HLambda}). Note that for any $i = 0, \dots, m$ and $T_i = T_{i,x_k}(\lambda_{i,k})$, we have
$$
\ba{rcl}
\lambda_{i,k} \la f'_i(x_k), x_k - T_i \ra & \refEQ{def-NPhi} & \vf_{i,x_k}(\lambda_{i,k}) + \beta_d(x_k,T_i) \; \refEQ{eq-HLambda} \; \half h_k^2 + \beta_d(x_k,T_i)\\
\\
& \refGE{eq-BGrow} & \half h_k^2 + \half \| x_k - T_i \|^2 \; \geq \; h_k \| x_k - T_i \|.
\ea
$$
Hence, since $T_i \neq x_k$,  we have
\beq\label{eq-LHLow}
\ba{rcl}
\lambda_{i,k} \| f'_{i}(x_k) \|_* & \geq & h_k, \quad k \geq 0.
\ea
\eeq

In method (\ref{met-SGM4}), we choose ${\cal H}$ in accordance to (\ref{def-HStep3}). Then, for $N \geq 1 + a(0)$, we define function $k(N)$ by the first equation in (\ref{def-KN}) and introduce by (\ref{def-Sigma}) the approximations of the optimal Lagrange multipliers.
\BT\label{th-Main3}
Let the sequence of points $\{ x_k \}_{k \geq 0}$ be generated by method (\ref{met-SGM4}). Then for any $N \geq 1 + a(0)$ the set ${\cal A}_0(N)$ is not empty. Moreover, if all functions $f_i(\cdot)$, $0 = 1, \dots, m$, satisfy (\ref{eq-GBound}), then for all $k \in {\cal A}_0(N)$ we have
\beq\label{eq-GSmal3}
\ba{rcl}
f_i(x_k) & \leq &  \| f'_{i}(x_k) \|_* h_k \; \refLE{eq-GBound} \; M_i h_{k(N)}, \quad i = 1 \dots m,
\ea
\eeq
\beq\label{eq-DGSmal3}
\ba{rcl}
{1 \over \sigma_0(N)} \sum\limits_{k \in {\cal A}_0(N)} \lambda_k f_0(x_k)  & \leq & \phi(\bar \lambda^*(N)) + M_0 h_{k(N)}.
\ea
\eeq

\ET
\proof
As in the proof of Theorem \ref{th-Main2}, for $x \in Q$, we denote $r_k(x) = \beta_d(x_k,x)$, and
$$
\ba{rcl}
B_N & \Def & \max\limits_{x \in Q} \Big[  \sum\limits_{k \in {\cal A}_0(N)} \lambda_k f_0(x_k)  - \sigma_0(N) f_0(x) - \sum\limits_{i=1}^p \sigma_i(N) f_i(x) \Big]\\
\\
& = & \max\limits_{x \in Q} \Big[  \sum\limits_{k \in {\cal A}_0(N)} \lambda_k [f_0(x_k)  - f_0(x)] - \sum\limits_{i=1}^m \sigma_i(N) f_i(x) \Big].
\ea
$$
Then, for  $k \in {\cal A}_0(N)$, we have proved that
$$
\ba{rcl}
r_{k+1}(x) - r_k(x) & \leq & - \lambda_k [f_0(x_k) - f_0(x)] + \half h_k^2.
\ea
$$

On the other hand,
$$
\ba{rcl}
\sum\limits_{i=1}^m \sigma_i(N) f_i(x) & \refEQ{def-Sigma} & \sum\limits_{i=1}^m \sum\limits_{k \in {\cal A}_i(N)} \lambda_k f_i(x) \; = \; 
\sum\limits_{k \not\in {\cal A}_0(N)} \lambda_{i_k,k} f_{i_k}(x)\\
\\
& \geq & \sum\limits_{k \not\in {\cal A}_0(N)} \lambda_{i_k,k} [f_{i_k}(x_k) + \la f'_{i_k}(x_k), x - x_k \ra ]\\
\\
& \stackrel{(\ref{met-SGM4})_2}{\geq} & 
\sum\limits_{k \not\in {\cal A}_0(N)} \Big[ h_k^2 + \lambda_{i_k,k} \la f'_{i_k}(x_k), x - x_k \ra \Big].
\ea
$$
Note that the first-order optimality condition for problem (\ref{def-Phi}) is as follows:
$$
\ba{rcl}
\la \nabla d(T_i) - \nabla d(x_k) + \lambda f'_i(x_k), x - T_i \ra & \geq & 0, \quad \forall x \in Q,
\ea
$$
where $T_i$ is its optimal solution, Therefore, for all $x \in Q$ an $k \not\in {\cal A}_0(N)$, we have
$$
\ba{rl}
& \lambda_{i_k,k} \la f'_{i_k}(x_k), x - x_k \ra \; = \; \lambda_{i_k,k} \Big[ \la f'_{i_k}(x_k), x - x_{k+1} \ra + \la f'_{i_k}(x_k), x_{k+1} - x_k \ra \Big]\\
\\
\geq & \lambda_{i_k,k} \la f'_{i_k}(x_k), x_{k+1} - x_k \ra + \la \nabla d(x_{k+1}) - \nabla d(x_k), x_{k+1} - x \ra\\
\\
= & \lambda_{i_k,k} \la f'_{i_k}(x_k), x_{k+1} - x_k \ra + r_{k+1}(x) - r_{k}(x) + \beta_d(x_k,x_{k+1})\\
\\
= & r_{k+1}(x) - r_{k}(x) - \vf_{i_k,x_k}(\lambda_{i,k}) \; \refEQ{eq-HLambda} \; r_{k+1}(x) - r_{k}(x) - \half h_k^2.
\ea
$$
Thus, as in the proof of Theorem \ref{th-Main2}, we conclude that
$$
\ba{rcl}
B_N & \leq & \max\limits_{x \in Q} \Big[ \sum\limits_{k \in {\cal A}_0(N)} [r_k(x) - r_{k+1}(x) + \half h_k^2] - \sum\limits_{k \not\in {\cal A}_0(N)} [ r_{k+1}(x) - r_{k}(x) + \half h_k^2 ] \; \Big]\\
\\
& \leq & h_{k(N)} \sum\limits_{k \in {\cal A}_0(N)} h_k \; \refLE{eq-LHLow} \; M_0 h_{k(N)} \sigma_0(N). 
\ea
$$
Since ${\cal A}_0(N) \neq \emptyset$ (see the proof of Theorem \ref{th-Main2}), dividing this inequality by $\sigma_0(N) > 0$, we get inequality (\ref{eq-DGSmal3}).

Finally, note that for all $k \in {\cal A}_0(N)$ and $i = 1, \dots, m$, we have 
$$
\ba{rcl}
\lambda_{i,k} f_i(x_k) & \leq & h_k^2 \; \refLE{eq-LHLow} \; \lambda_{i,k} \| f'_i(x_k) \|_* h_k \; \refLE{eq-GBound} \; \lambda_{i,k} M_i h_k \; \leq \; \lambda_{i,k} M_i h_{k(N)} .
\ea
$$
Thus, inequality (\ref{eq-GSmal3}) is proved.
\qed

\section{Accuracy guarantees for the dual problem}\label{sc-AcDual}
\SetEQ

In Sections \ref{sc-Lagrange1} and \ref{sc-Lagrange2}, we developed two convergent methods (\ref{met-SGM3}) and (\ref{met-SGM4}), which are able to approach the optimal solution of the primal problem (\ref{prob-FConst}), generating in parallel an approximate solution of the dual problem (\ref{prob-LDual}). Indeed, for $N \geq 1 + a(0)$, denote
$$
\ba{rcl}
f_0^*(N) & = & \min\limits_{k \in {\cal A}_0(N)} f_0(x_k), \quad x^*_N \; = \; \arg\min\limits_{x} \Big\{ f_0(x): \; x = x_k, \; k \in {\cal A}_0(N) \Big\}.
\ea
$$
Then, in view of Theorems \ref{th-Main2} and (\ref{th-Main3}), we have
\beq\label{eq-BGap}
\ba{rcl}
f_0^*(N) - \phi(\bar \lambda^*(N)) & \leq & M_0 h_{k(N)},\\
\\
\max\limits_{1 \leq i \leq m} {1 \over M_i} f_i(x_k) & \leq & h_{k(N)}, \quad k \in {\cal A}_0(N).
\ea
\eeq
Since $\phi(\bar \lambda^*(N)) \leq f^*_0$, this inequality justfies that the point $x^*_N$ is a good approximate solution to primal problem (\ref{prob-FConst}). 

However, note that in our reasonings, we did not assume yet the {\em existence} of optimal solution of the dual problem (\ref{prob-LDual}). It appears that under our assumptions, this may not happen. In this case, since $f_0^*(N)$ can be {\em significantly smaller} than $f_0^*$, inequality~(\ref{eq-BGap}) cannot justify that vector $\bar \lambda^*(N)$ delivers a good value of the dual objective function.

Let us look at the following example.
\BE\label{ex-NoSlater}
Consider the following problem 
$$
\ba{c}
\min\limits_{x \in \R^2}
\Big\{ x^{(2)}:\; 1 -x^{(1)}\leq 0, \; \| x \|_2 \leq 1 \Big\}.
\ea
$$
In this case, ${\cal L}(x,\lambda) = x^{(2)} + \lambda (1
- x^{(1)})$. Hence,
$$
\ba{rcl}
\phi(\lambda) & = & \min\limits_{\| x \|_2 \leq 1} {\cal
L}(x,\lambda) \; = \; \lambda - [1 + \lambda^2]^{1/2}.
\ea
$$
Thus, there is no duality gap: $\phi_* \Def \sup\limits_{\lambda \geq
0} \phi(\lambda) = 0$. However, the optimal dual solution
$\lambda^*$ does not exist.

Let us look now at the perturbed feasible set
$$
\ba{rcl}
{\cal F}()\epsilon) & = & \{ x \in \R^2: \; \| x \|_2 \leq 1, \; 1 - \epsilon \leq x^{(1)} \},
\ea
$$ 
where $\epsilon > 0$ is sufficiently small. Note that it contains a point with the second coordinate equal to $\phi_* - \sqrt{\epsilon(2-\epsilon)}$. This means that the condition (\ref{eq-BGap}) can guarantee only that
$$
\ba{rcl}
\phi(\lambda^*(N) & \geq & \phi_* - \epsilon - \sqrt{\epsilon(2-\epsilon)}.
\ea
$$
Hence, for dual problems with nonexisting optimal solutions, we can expect a significant drop in the quality of approximation in terms of the function value.
\qed
\EE

Thus, in our complexity bounds, we need to take into account the size of  optimal dual solution $\bar \lambda^* \in \R^m_+$. Let the sequence $\{ x_k \}_{k \geq 0}$ be generated by one of the methods (\ref{met-SGM3}) or~(\ref{met-SGM4}). Then, for any $k \geq 0$, we have
\beq\label{eq-FDLow}
\ba{rl}
f_0(x_k) \geq & \min\limits_{x \in Q} \Big\{ f_0(x): \; \bar f(x) \leq \bar f(x_k) \Big\} = \min\limits_{x \in Q} \max\limits_{\bar \lambda \in \R^m_+} \Big\{ f_0(x) + \la \bar \lambda, \bar f(x) - \bar f(x_k) \ra \Big\}\\
\\
= & \max\limits_{\bar \lambda \in \R^m_+}  \min\limits_{x \in Q} \Big\{ f_0(x) + \la \bar \lambda, \bar f(x) - \bar f(x_k) \ra \Big\} \; = \; \max\limits_{\bar \lambda \in \R^m_+}  \Big\{ \phi(\bar \lambda) - \la \bar \lambda, \bar f(x_k) \ra \Big\}\\
\\
\geq & f^*_0 - \la \bar \lambda^*,  \bar f(x_k) \ra \; \geq \; f^*_0 - \la \bar \lambda^*, \bar M \ra \max\limits_{1 \leq i \leq m} {1 \over M_i} f_i(x_k).
\ea
\eeq
Thus, we have proved the following theorem.
\BT\label{th-QDual}
Under conditions of Theorem \ref{th-Main2} or \ref{th-Main3}, for all $N \geq 1 + a(0)$ and all $k \in {\cal A}_0(N)$, we have
\beq\label{eq-DRate}
\ba{rcl}
f_0^* - \phi(\bar \lambda^*(N)) & \leq & (M_0 + \la \bar \lambda^*, \bar M \ra) \, h_{k(N)}.
\ea
\eeq
\ET
\proof
Indeed, in view of inequalities (\ref{eq-GSmall}) and (\ref{eq-GSmal3}), we have $\max\limits_{1 \leq i \leq m} {1 \over M_i} f_i(x_k) \leq h_{k(N)}$. Thus, we get the bound (\ref{eq-DRate}) from (\ref{eq-BGap}) and (\ref{eq-FDLow}).
\qed

Recall that the size of optimal dual multipliers can be bounded by the standard {\em Slater condition}. Namely, let us assume existence of a point $\hat x \in Q$ such that
\beq\label{eq-CSlater}
\ba{rcl}
f_i(\hat x) & < & 0, \quad i = 1, \dots, m.
\ea
\eeq
Then, in accordance to Lemma 3.1.21 in \cite{NewLec}, we have
\beq\label{eq-MBound}
\ba{rcl}
\la \bar \lambda^*, - \bar f(\hat x) \ra & \leq & f_0(\hat x) - f^*_0.
\ea
\eeq
Therefore,
$$
\ba{rcl}
\la \bar \lambda^*, \Bar M \ra & = &
\sum\limits_{i=1}^m \lambda^*_i(-f_i(\hat x)) \cdot {M_i \over - f_i(\hat x)} \; \leq \; \la \bar \lambda^*, - \bar f(\hat x) \ra \max\limits_{1 \leq i \leq m} {M_i \over - f_i(\hat x)}\\
\\
& \refLE{eq-MBound} & (f_0(\hat x) - f^*_0) \max\limits_{1 \leq i \leq m} {M_i \over - f_i(\hat x)}.
\ea
$$
Thus, the Slater condition provides us with the following bound:
\beq\label{eq-BSlater}
\ba{rcl}
f_0^* - \phi(\bar \lambda^*(N)) & \leq & \Big(M_0 + (f_0(\hat x) - f^*_0) \max\limits_{1 \leq i \leq m} {M_i \over - f_i(\hat x)} \Big) \, h_{k(N)},
\ea
\eeq
which is valid for all $N \geq 1 + a(0)$. Note that we are able to compute vector $\bar \lambda^*(N)$ without computing values of the dual function $\phi(\cdot)$, which can be very complex. In fact, computational complexity of a single value $\phi(\lambda)$ can be of the same order as the complexity of solving the initial problem (\ref{prob-FConst}), or even more.

\section{Subgradient method for unbounded feasible set}\label{sc-NoBound}
\SetEQ

In the previous sections, we looked at optimization methods applicable to the  bounded sets (see condition (\ref{eq-DHBound})). If this is not true, the second-order divergence condition (\ref{eq-TGrow}) cannot help, and we need to find another way of justifying efficiency of the subgradient schemes. This is the goal of the current section.

We are still working with the problem (\ref{prob-FConst}), satisfying condition (\ref{eq-GBound}). However, the set $Q$ is not bounded anymore. 
Hence, we cannot count on Sion's theorem (\ref{prob-LDual}).

In our method, we have a sequence of scaling coefficients $\Gamma = \{ \gamma_k \}_{k \geq 0}$ satisfying condition
\beq\label{eq-CondG}
\ba{rcl}
\gamma_{k+1} & > & \gamma_k \; \geq \; 0, \quad k \geq 0,
\ea
\eeq
a tolerance parameter $\epsilon > 0$, and a rough estimate for the distance to the optimum $D_0 > 0$. 

For functions $ \vf_{i,y}(\cdot)$ defined by (\ref{def-NPhi}) with $y \in Q$, denote by $a = a_{i,k}(y)$ the unique solution of the equation
\beq\label{eq-NStep}
\ba{rcl}
\vf_{i,y}\left( {a\over \gamma_{k+1}} \right) & = & \left(1 - {\gamma_k \over \gamma_{k+1}} \right) D_0, \quad i = 0, \dots, m.
\ea
\eeq

Let us look at the following optimization scheme.
\beq\label{met-SGM5}
\ba{|c|}
\hline \\
\quad \mbox{\bf Switching Subgradient Method for Unbounded Sets} \quad\\
\\
\hline \\
\ba{l}
\mbox{{\bf Initialization.}  Choose $x_0 \in Q$, a sequence $\Gamma$ (see (\ref{eq-CondG})), and some $D_0 > 0$.}\\
\\
\mbox{{\bf $k$th iteration ($k \geq 0$).} {\bf 1.} Define $\tau_k = {\gamma_k \over \gamma_{k+1}}$ and $y_k = (1-\tau_k)x_0 + \tau_k x_k$.} \\
\\
\mbox{{\bf 2.} If $\exists \; i \in \{ 1, \dots, m\}:\; f_i(y_k) \geq \epsilon$, then set $i_k = i$. Otherwise, set $i_k = 0$.}\\
\\
\mbox{{\bf 3.} Compute $a_k = a_{i_k,k}(y_k)$ by (\ref{eq-NStep}), and set $x_{k+1} = T_{i_k,y_k}\Big( {a_k \over \gamma_{k+1}} \Big)$.}\\
\\
\ea\\
\hline
\ea
\eeq
It seems that now, the selection of violated constraint by Step 2 looks more natural than in the methods (\ref{met-SGM3}) and (\ref{met-SGM4}).

For $x \in Q$, denote $\Delta_k(x) = \gamma_k\Big(\beta_d(x_k,x) - \beta_d(x_0,x)\Big)$. This is a linear function of~$x$. In what follows, we assume the Bregman distance is convex with respect to its first argument:
\beq\label{eq-BConv}
\ba{rcl}
\beta_d(\alpha u_1 + (1-\alpha)u_2,x) & \leq & \alpha \beta_d(u_1,x) + (1- \alpha) \beta_d(u_2,x), \\
\\
& & \forall u_1,u_2, x \in Q, \; \alpha \in [0,1].
\ea
\eeq
This property is not very common. However, it is valid for the following two important examples.
\BE\label{ex-BConv}$  $

\noindent
{\bf 1.} Let $d(x) = \half \| x \|^2_B$ ( see (\ref{def-Euclid})). Then $\beta_d(x,y) = \half \| x - y \|^2_B$ and (\ref{eq-BConv}) holds.

\noindent
{\bf 2.} Let $Q = \{ x \in \R^n_+: \; \sum\limits_{i=1}^n x^{(i)} = 1 \}$. Define $d(x) = \sum\limits_{i=1}^n x^{(i)} \ln x^{(i)}$. Then, for $x, y \in Q$, we have $\beta_d(x,y) = \sum\limits_{i=1}^n y^{(i)} \ln {y^{(i)} \over x^{(i)}}$, and (\ref{eq-BConv}) holds also.
\qed
\EE

Now we can prove the following statement.
\BL\label{lm-Change}
Let the Bergman distance $\beta_d(\cdot,\cdot)$ satisfy (\ref{eq-BConv}). Then, for the sequence $\{ x_k \}_{k \geq 0}$ generated by method (\ref{met-SGM5}), and any  $x \in Q$, we have
\beq\label{eq-Change}
\ba{rcl}
\Delta_{k+1}(x) & \leq & \Delta_k(x) + a_k \la f'_{i_k}(y_k), x - y_k \ra + (\gamma_{k+1} - \gamma_k) D_0, \quad k \geq 0.
\ea
\eeq
\EL
\proof
Note that
$$
\ba{rcl}
\Delta_{k+1}(x) - \Delta_k(x) & = & \gamma_{k+1} \Big( \beta_d(x_{k+1},x) - \beta_d(x_0,x) \Big) - \gamma_{k} \Big( \beta_d(x_{k},x) - \beta_d(x_0,x) \Big)\\
\\
& \refLE{eq-BConv} & \gamma_{k+1} \Big( \beta_d(x_{k+1},x) - \beta_d(y_k,x) \Big)\\
\\
& \refEQ{def-Breg} & \gamma_{k+1} \Big( \la \nabla d(y_k) - \nabla d(x_{k+1}), x - x_{k+1} \ra - \beta_d(y_k,x_{k+1} \ra \Big).
\ea
$$
Since the first-order optimality conditions for point $x_{k+1}$ tells us that
$$
\ba{rcl}
\la a_k  f'_{i_k}(y_k) + \gamma_{k+1}(\nabla d(x_{k+1}) - \nabla d(y_k)), x - x_{k+1} \ra & \geq & 0, \quad x \in Q,
\ea
$$
we have
$$
\ba{rcl}
\Delta_{k+1}(x) - \Delta_k(x) & \leq & a_k \la f'_{i_k}(y_k), x - x_{k+1} \ra - \gamma_{k+1} \beta_d(y_k, x_{k+1} \ra\\
\\
& \refEQ{def-NPhi} & a_k \la f'_{i_k}(y_k), x - y_k \ra + \gamma_{k+1} \vf_{i,y_k} \left( {a_k \over \gamma_{k+1}}\right) \\
\\
& \refEQ{eq-NStep} & a_k \la f'_{i_k}(y_k), x - y_k \ra + (\gamma_{k+1} - \gamma_k) D_0,
\ea
$$
and this is inequality (\ref{eq-Change}).
\qed

Since $\Delta_0(x) = 0$, we can sum up the inequalities (\ref{eq-Change}) for $k = 0, \dots, N-1$ with $N \geq 1$, and get the following consequence:
\beq\label{eq-BSum}
\ba{rcl}
\gamma_N \beta_d(x_N, x) + \sum\limits_{k=0}^{N-1} a_k \la f'_{i_k}(y_k), y_k - x \ra & \refLE{eq-Change} & \gamma_N ( \beta_d(x_0,x) + D_0), 
\ea
\eeq
which is valid for all $x \in Q$. 

In order to approximate the optimal Lagrange multipliers of problem (\ref{prob-LDual}), we need to introduce the following objects:
\beq\label{def-LMult}
\ba{rcl}
{\cal S}_i(N) & = & \{ k: \; 0 \leq k \leq N-1, \; i_k = i \}, \quad N \geq 1,\\
\\
\hat \sigma_i(N) & = & \sum\limits_{k \in {\cal S}_i(N)} a_k, \quad i = 0, \dots, m,\\
\\
\hat \lambda_i(N) & = & \hat \sigma_i(N)/ \hat \sigma_0(N), \quad i = 1, \dots, m.
\ea
\eeq
Denote $\hat \lambda(N) = (\hat \lambda_1(N), \dots, \hat \lambda_m(N))^T \in \R^m_+$ and $f_0^*(N) = \min\limits_{k \in {\cal S}_0(N)} f_0(x_k)$. Note that for all $k \in {\cal S}_0(N)$ we have
\beq\label{eq-FFeas}
\ba{rcl}
f_i(x_k) & \leq & \epsilon, \quad k = 1, \dots, m.
\ea
\eeq

For our convergence result, we need to assume existence of the optimal solution $x^*$ to problem (\ref{prob-FConst}). Let us introduce some bound
$$
\ba{rcl}
D & \geq & \beta_d(x_0,x^*),
\ea
$$
which is not used in the method(\ref{met-SGM5}).
Denote $Q_D = \{ x \in Q: \; \beta_d(x_0,x) \leq D \}$. Clearly, if we replace in problem (\ref{prob-FConst}) the set $Q$ by $Q_D$, then its optimal solution will not be changed. However, now we can correctly
define a {\em restricted} dual function $\phi_D(\lambda) = \min\limits_{x \in Q_D} {\cal L}(x,\lambda)$. Let us prove our main result.

\BT\label{th-Main4}
Let functional components of problem (\ref{prob-FConst}) satisfy the following condition:
\beq\label{eq-UGBound}
\ba{rcl}
\| f'_i(x) \|_* & \leq & M, \quad i = 1, \dots, m.
\ea
\eeq
Then, as far as
\beq\label{eq-IBound}
\ba{rcl}
\Sigma_N \; \Def \; {1 \over \gamma_N} \sum\limits_{k=0}^{N-1} \Big[\gamma_{k+1}(\gamma_{k+1} - \gamma_k)\Big]^{1/2} & > &  {\rho_0 + D_0 \over  \sqrt{2D_0}} \cdot {M \over \epsilon},
\ea
\eeq
where $\rho_0 = \min\limits_{y \in {\cal F}}\beta_d(x_0,y)$,
we have ${\cal S}_0(N) \neq \emptyset$ and $\hat \sigma_0(N) > 0$. Moreover, if
\beq\label{eq-NBig}
\ba{rcl}
\Sigma_N & > & {D + D_0 \over  \sqrt{2D_0}} \cdot {M \over \epsilon},
\ea
\eeq
then
\beq\label{eq-SGap}
\ba{rcl}
f^*_0(N) - f^*_0 \; \leq \; f^*_0(N) -  \phi_D(\hat \lambda(N)) & \leq & \epsilon.
\ea
\eeq
\ET
\proof
Indeed, for any $x \in Q$, we have
$$
\ba{rcl}
\sum\limits_{k=0}^{N-1} a_k \la f'_{i_k}(y_k), y_k - x \ra & = & \sum\limits_{k \in {\cal S}_0(N)} a_k \la f'_0(y_k), y_k - x \ra + \sum\limits_{i=1}^m \sum\limits_{k \in {\cal S}_i(N)} a_k \la f'_i(y_k), y_k - x \ra\\
\\
& \geq & \sum\limits_{k \in {\cal S}_0(N)} a_k [f_0(y_k) - f_0(x)] 
+ \sum\limits_{i=1}^m \sum\limits_{k \in {\cal S}_i(N)} a_k [f_i(y_k) - f_i(x)]\\
\\
& \geq & \hat \sigma_0(N) [f^*_0(N) - f_0(x)] + \sum\limits_{i=1}^m \hat \sigma_i(N)[\epsilon - f_i(x)].
\ea
$$
Hence, in view of inequality (\ref{eq-BSum}), we get the following bound
$$
\ba{rcl}
\hat \sigma_0(N) [f^*_0(N) - f_0(x)] + \sum\limits_{i=1}^m \hat \sigma_i(N)[\epsilon - f_i(x)] & \leq & \gamma_N(\beta_d(x_0,x) + D_0).
\ea
$$
Note that $\hat\sigma_0(N) + \sum\limits_{i=1}^m \hat \sigma_i(N) = \sum\limits_{k=0}^{N-1} a_k$. Therefore, this inequality can be rewritten as follows:
\beq\label{eq-IMain}
\ba{rcl}
\gamma_N(\beta_d(x_0,x) + D_0) & \geq & \hat \sigma_0(N) [ f^*_0(N)- f_0(x) - \epsilon] + \epsilon \sum\limits_{k=0}^{N-1} a_k - \sum\limits_{i=1}^m \hat \sigma_i(N) f_i(x).
\ea
\eeq
In view of condition(\ref{eq-NStep}), for $T_i = T_{i,y_k}\left( {a_{i,k}(y_k) \over \gamma_{k+1}} \right)$, we have
$$
\ba{rl}
& a_{i,k}(y_k) \la f'_i(y_k), y_k - T_i \ra \; = \; \gamma_{k+1} \beta_d(y_k,T_i)  +  (\gamma_{k+1} - \gamma_k) D_0\\
\\
\refGE{eq-BGrow} & \half \gamma_{k+1} \| T_i - y_k \|^2 + (\gamma_{k+1} - \gamma_k) D_0 \; \geq \; \Big[ 2 \gamma_{k+1} (\gamma_{k+1} - \gamma_k) D_0 \Big]^{1/2} \| T_i - y_k \|.
\ea
$$
Thus, $M a_k(y_k) \refGE{eq-UGBound} \Big[ 2 \gamma_{k+1} (\gamma_{k+1} - \gamma_k) D_0 \Big]^{1/2}$, and we conclude that 
\beq\label{eq-SigmaN}
\ba{rcl}
\sum\limits_{k=0}^{N-1} a_k & \geq & {1 \over M}\sqrt{2D_0} \gamma_N \Sigma_N.
\ea
\eeq

Let us assume that ${\cal S}_0(N) = \emptyset$. Then $\hat \sigma_0(N) = 0$ and for $x  = \arg\min\limits_{y \in {\cal F}} \beta_d(x_0,y)$ inequality (\ref{eq-IMain}) leads to the following relation:
$$
\ba{rcl}
\gamma_N( \rho_0 + D_0) & \refGE{eq-SigmaN} & {\epsilon \over M} \sqrt{2D_0} \gamma_N \Sigma_N.
\ea
$$
However, this cannot happen in view of condition (\ref{eq-IBound}). Hence, inequality (\ref{eq-IBound}) implies ${\cal S}_0(N) \neq \emptyset$ and $\hat \sigma_0(N) > 0$. In this case, inequality (\ref{eq-IMain}) can be rewritten as follows:
$$
\ba{rcl}
\beta_d(x_0,x) + D_0 & \refGE{eq-SigmaN} & {1 \over \gamma_N}\hat \sigma_0(N) [ f^*_0(N)- {\cal L}(x, \hat \lambda(N)) - \epsilon] + {\epsilon \over M}\sqrt{2D_0} \Sigma_N.
\ea
$$
Maximizing the right-hand side of this inequality in $x \in Q_d$, we get
$$
\ba{rcl}
D + D_0 & \geq & {1 \over \gamma_N}\hat \sigma_0(N) [ f^*_0(N)- \phi_D(\hat \lambda(N)) - \epsilon] + {\epsilon \over M}\sqrt{2D_0} \Sigma_N.
\ea
$$
Hence, inequality (\ref{eq-NBig}) implies (\ref{eq-SGap}).
\qed

Thus, for the fast rate of convergence of method (\ref{met-SGM5}), we need to ensure a fast growth of the valies $\Sigma_N$. Let us choose
\beq\label{def-Gamma}
\ba{rcl}
\gamma_k & = & \sqrt{k}, \quad k \geq 0.
\ea
\eeq 
Then $\gamma_{k+1}(\gamma_{k+1} - \gamma_k) = {\sqrt{k+1} \over \sqrt{k+1} + \sqrt{k}} \geq \half$. Hence,
\beq\label{eq-LSigma}
\ba{rcl}
\Sigma_N & \geq & {1 \over \sqrt{N}} \cdot {N \over \sqrt{2}} \; = \; \sqrt{N \over 2}.
\ea
\eeq
In this case, inequality (\ref{eq-SGap}) can be ensured in $O(\epsilon^{-2})$ iterations.

Note that method (\ref{met-SGM5}) is quite different from the existing optimization schemes. Let us write down how it looks like in the case $Q = \E$, $d(x) = \half \| x \|^2_B$, and the parameter choice (\ref{def-Gamma}). In this case, the equation (\ref{eq-NStep}) can be written as follows:
$$
\ba{rcl}
{a^2 \over 2 \gamma_{k+1}^2} \| f'_i(y_k) \|^2 & = & {\gamma_{k+1} - \gamma_k \over \gamma_{k+1}} D_0.
\ea
$$
Thus, $a_k = \sqrt{2 \gamma_{k+1}(\gamma_{k+1} - \gamma_k) D_0 } \cdot {1 \over \| f'_{i_k}(y_k) \|^*_B} \approx {\sqrt{D_0} \over \| f'_{i_k}(y_k) \|^*_B}$, and the method looks as follows:
\beq\label{met-SGM6}
\ba{|rl|}
\hline & \\
& \hspace{8ex} \mbox{\bf Switching Subgradient Method With $Q = \E$}\\
& \\
\hline & \\
\mbox{\bf 0.} & \mbox{Choose $x_0 \in \E$ and $D_0 > 0$. {\bf For $k \geq 0$, iterate:}}\\
& \\
\mbox{\bf 1.} & \mbox{Define $\tau_k = \sqrt{k \over k+1}$ and $y_k = (1-\tau_k) x_0 + \tau_k x_k$.}\\
& \\
\mbox{\bf 2.} & \mbox{If $\exists i: 1 \leq i \leq m, \; f_i(y_k) \geq \epsilon$, then set $i_k = i$. Otherwise, set $i_k = 0$.}\\
& \\
\mbox{\bf 3.} & \mbox{Set $x_{k+1} = x_k - \sqrt{D_0 \over k+1} \cdot B^{-1} f'_{i_k}(y_k)/ \| f'_{i_k} (y_k) \|^*_B$.}\\
& \\
\hline
\ea
\eeq

\end{document}